\documentclass[10pt]{article}
\usepackage{amssymb, amsmath, latexsym, a4}
\usepackage[latin1]{inputenc}
\usepackage{graphicx}
\usepackage[all]{xy}
\usepackage[usenames]{color}
\usepackage{pgf,tikz}

\usepackage{mathrsfs}

\usetikzlibrary{arrows}

\newenvironment{proof}[1][Proof]{\noindent\textbf{#1.} }{\ \rule{0.5em}{0.5em}}

\newtheorem{De}{Definition}[section]
\newtheorem{Th}[De]{Theorem}
\newtheorem{Pro}[De]{Proposition}

\newtheorem{Rem}[De]{Remark}
\newtheorem{Ex}[De]{Example}

\newbox\pullbackbox
\setbox\pullbackbox=\hbox{\xy 0;<1mm,0mm>: \POS(4,0)\ar@{-} (0,0)
\ar@{-} (4,4)
\endxy}

\newdir{>>}{{}*!/3.5pt/:(1,-.2)@^{>}*!/3.5pt/:(1,+.2)@_{>}*!/7pt/:(1,-.2)@^{>}*!/7pt/:(1,+.2)@_{>}}
\newdir{ >>}{{}*!/8pt/@{|}*!/3.5pt/:(1,-.2)@^{>}*!/3.5pt/:(1,+.2)@_{>}}
\newdir{ |>}{{}*!/-3.5pt/@{|}*!/-8pt/:(1,-.2)@^{>}*!/-8pt/:(1,+.2)@_{>}}
\newdir{ >}{{}*!/-8pt/@{>}}
\newdir{>}{{}*:(1,-.2)@^{>}*:(1,+.2)@_{>}}
\newdir{<}{{}*:(1,+.2)@^{<}*:(1,-.2)@_{<}}

\begin{document}
\centerline{\bf Simplicial Structure on Connected Multiplicative
Operads}

\bigskip
\centerline{  Batkam Mbatchou V. Jacky III$^{1};$  Calvin
Tcheka$^{2}.$}

\bigskip

\centerline{$^{1}$ Department of Mathematics, Faculty of
    science-University of Dschang} \centerline{Campus Box (237)67
    Dschang, Cameroon} \centerline{ {E-mail address}: batkamjacky3@yahoo.com }

\bigskip

\centerline{$^{2}$ Department of Mathematics, Faculty of
science-University of Dschang} \centerline{Campus Box (237)67
Dschang, Cameroon} \centerline{ {E-mail address}:
calvin.tcheka@univ-dschang.org}

\bigskip
\date{}

\bigskip \bigskip \bigskip

 \centerline{\bf Abstract}  In these notes, we define a new simplicial
 structure on a connected multiplicative operad and call it {\it  connected
 multiplicative simplicial operad (for short; simplicial operad)}.
 Next we introduce on this simplicial operad a brace algebra structure analogous to that of
Gerstenhaber-Voronov that we call {\it right brace algebra
structure.} This permits us to obtain on the operad with the above
mentioned properties a bicomplex structure one of whose two
differential operators is a coboundary and the other one is a
boundary. Moreover we define on one hand on the above simplicial
operad together with its right brace algebra structure, two distinct
products up to a sign respectively called {\it dot-product} and {\it
odot-product}. Then we show that the coboundary and the boundary
together with the {\it odot-product} provide to this simplicial
operad two distinct differential graded algebra structures. On the
other hand we obtain through the Alexander-Withney
 map, a differential graded coalgebra structure on a simplicial operad.
We end by illustrating our constructions with some examples.

\bigskip

{\bf Keywords:}  { Simplicial operad, Simplicial right brace
algebra, Simplicial coalgebra, dot-product, odot-product, Hochschild
cohomology, Shift operad.}

\bigskip

 {\bf MSC:} 16N60, 18N50, 18G31, 19D23, 18M05, 57T99.
\bigskip

\begin{center}\section{Introduction}\end{center}
Operads are tools usually used to describe algebraic structures in
monoidal categories. They are very important in categories with a
good notion of homotopy where they are useful for the study of
homotopy invariant algebraic structures and hierarchies of higher
homotopies. Although one can already find its trace in the paper of
Lazard \cite{La} entitled group laws and analyzers, the basic idea
of operads has mainly been developed in Chicago in the seventies by
the algebraic topologists(S. MacLane \cite{McL}, J. Stasheff
\cite{St}, J.P. May \cite{May}, J.M. Boardmann, R. Vogt \cite{B.
V.}, F.R. Cohen \cite{Coh}) to study loop spaces. In addition
instead of describing algebras by its generators(operations) and the
relatives relations(fundamental identities), one can consider all
operations that can be performed on a finite number of variables and
the relations between these operations. This structure has been
baptized by J. P. May: Operad structure. The main interest of this
point of view resides on the fact that one can compare algebras even
if they are of different natures. In other words one has the notion
of homomorphisms between operads.\\  In the litterature the
expression simplicial operad denotes an operad in the category of
simplicial sets, {\it sSet}(see \cite{1}, \cite{6}). In this work we
consider connected multiplicative operads in the category of
$\mathbb{K}$-vector spaces and introduce a new variant of simplicial
operad structure. Indeed we define the notions of face map and
degeneracy map which lead us to the construction of the above
mentioned simplicial operad structure. Next we introduce on a
connected multiplicative operad a brace algebra structure analogous
to that of Gerstenhaber-Voronov and call it {\it right brace algebra
structure.} We specify that our brace algebra structure differs from
the one previously constructed on how to count the inputs around a
given operation in the brace. Subsequently, we define by using this
new brace algebra structure, a product called dot-product analogous
to that of Gerstenhaber-Voronov(see 4.2) which later will be
slightly modified to obtain the one called {\it odot-product} and
denoted $\odot$(see once more 4.2). Afterward, we define through the
face map on this simplicial connected multiplicative operad a
derivation, $\partial,$ with respect to the modified dot-product
such that $((\mathcal{O}, \partial), \odot)$ becomes a differential
graded algebra. The motivation for the definition of our right brace
structure is based on the following fact: given a unitary
associative $\mathbb{K}$-algebra A, we want our constructions to be
conform to some well known result of the sub operad
$\mathcal{L}_{A}^{\prime}$ of the endomorphism operad, ${\it
End}_{A,}$ which encodes the Hochschild cochain complex of A, $\mathfrak{C}^{\ast}(A, A)$(see 5.1).\\
 Beside the above mentioned constructions, we define a variant of Alexander-Withney
 map on a given simplicial operad $\mathcal{O}=\{O(n), n\ge 0\}$ by using the face map.
 This leads us to the construction of a coproduct, $\triangle_{\mathcal{O}},$ on $\mathcal{O}$  such that the triple
 $((\mathcal{O}, \partial), \Delta_{\mathcal{O}}, \epsilon)$ is a co-unitary differential graded
coalgebra(see subsection 3.2).\\ For the sake of illustration of the
above mentioned constructions, we focused our attention on two
particular operads, namely the endomorphism operad, $\mbox{End}_{A}$
and the associative operad, $\mathcal{A}_{ss}.$ Beside these two
operads,
we also define a new operad called shift operad on which we construct a simplicial structure.\\
Firstly, it is well known that for a given ground field
$\mathbb{K},$ an associative unitary $\mathbb{K}$-algebra $A$ gives
rise to the endomorphism operad also denoted $\mathcal{L}_{A}:=
\bigoplus_{n\geq o}\mathcal{L}_{A}(n)= \bigoplus_{n\geq o}
Hom_{\mathbb{K}}(A^{\otimes n}, A)$ which is a non connected unitary
multiplicative $\Sigma$-operad and encodes the Hochschild cochain
complex, $\mathfrak{C}(A;A),$ of $A$ with in itself. However there
is a sub operad $\mathcal{L}^{\prime}_{A}:=
(\mathcal{L}^{\prime}_{A}(n))_{n\geq 0}$ of the operad
$\mathcal{L}_{A}$ such that
$\mathcal{L}^{\prime}_{A}(n)=\left\{\begin{array}{ll}
    \mathcal{L}^{\prime}_{A}(n)\quad \mbox{for all} \quad n\geq 1, \\
   \mathbb{K}\qquad \mbox{if} \quad n= 0
\end{array}
\right.$ which is a connected unitary multiplicative
$\Sigma$-operad, thus it is endowed with a simplicial right brace
algebra structure. Moreover the dot-product induced by its
associated right brace structure and the differential operator $d$
coincide respectively with the
ordinary cup-product and the Hochschild coboundary operator.\\
Secondly, we construct the new simplicial structure on the
associative operad $\mathcal{A}_{ss}.$ More precisely, if we denote
by $\sum= \bigoplus_{n \geq 1}\sum_{n},$ the collection of
permutations with $\sum_{n},$ the set of permutations of order $n\in
\mathbb{N}^{\ast};$ we define a partial composition of permutations
on $\sum.$ Using the connectivity and the multiplicative structure
of the operad $\mathcal{A}_{ss},$ this partial composition permits
us to define face maps and  degeneracy maps which together with the
brace algebra structure induce a differential graded algebra
structure on $\mathcal{A}_{ss}.$ Furthermore we show that the
coalgebra structure on the simplicial operad $\mathcal{A}_{ss}$
given by the above mentioned coproduct
$\triangle_{\mathcal{A}_{ss}}$
coincides with the Malvenuto-Reutenauer coalgebra structure(see \cite{C.P}, page 217 or \cite{M.R}).         \\
In the last example, we consider the set
$\mathcal{E}=\{E_{n}\}_{n\in \mathbb{N}}$ such that for $n\geq 1,$
$E_n=\{X_{n} = \{a_1\le\cdots\le a_n\} \mid a_i\in
\mathbb{N}^{\ast}, a_i\neq a_{i+1}, 1\leq i\leq n\}$ is a set of sub
sets of n distinct ordered elements of $ \mathbb{N}^{\ast}$ with
respect to the relation $\leq$ with the convention $E_0=\{X_{0} = 0=
\emptyset\}.$ We define on $\mathcal{E}$ an external operation
called shift operation through
which we build a simplicial structure.\\
 This work is organized as follows: in section
2, we give general recollection on operads. In section 3, we define
our variant of simplicial operad structure as well as the
derivations which lead us to a differential graded chain algebra. We
end this section with the construction of a coalgebra structure on
the simplicial operad through the Alexander-Whitney map. In section
4 we endow the simplicial operad with the right brace algebra
structure. Section 5
provides some examples for the illustration of our constructions.\\
\vspace{2cm}
\section{General reminders on operads}
Let $\mathbb{K}$ be an arbitrary ground field. Throughout these
notes, Operads considered are over the monoidal category of
$\mathbb{K}$-vector spaces denoted $\mbox{Vect}_{\mathbb{K}}.$ Such
operads are said to be symmetric($\Sigma$-operad) if they are
endowed with a right action of symmetric groups $\Sigma=
\{\Sigma_{n}, n\in \mathbb{N}\} $ and non symmetric(non
$\Sigma$-operad) if not(see \cite{Lo-Va}). In this section, we
recall the fundamental notions on operads over the category
$\mbox{Vect}_{\mathbb{K}}$ and their related properties. For more
details reader can refer to [\cite{B. D.}, \cite{F.}, \cite{10},
\cite{Lo-Va} ].
\begin{De}
{\it A symmetric operad($\Sigma$-operad) is a collection of right
$\sum$-vector spaces over $\mathbb{K},$ $\{\mathcal{O}(k)| k\geq
1\},$ together with a composition product:
$$\begin{array}{ccc}
 \mathcal{O}(k) \otimes \mathcal{O}(n_{1}) \otimes \cdots \otimes \mathcal{O}(n_{k})& \stackrel{\gamma^{\mathcal{O}}}\longrightarrow & \mathcal{O}(n_{1} +\cdots + n_{k})\\
x \otimes x_1 \otimes \cdots \otimes x_k & \longmapsto &
\gamma^{\mathcal{O}}(x; x_1,\cdots , x_k)
\end{array}$$
which is:
\begin{enumerate}
\item associative in the sense that \\
$ \begin{array}{ccc} \gamma^{\mathcal{O}}(\gamma^{\mathcal{O}}(x;
x_1, \cdots , x_k); y_1,\cdots, y_{n_1+ \cdots +n_k}) & =&
\gamma^{\mathcal{O}}(x; \gamma^{\mathcal{O}}(x_1; y_1,\cdots ,
y_{n_1}), \gamma^{\mathcal{O}}(x_2;  y_{n_{1}+ 1},\cdots , y_{n_{1}+
n_{2}})\\, &  &. . . , \gamma^{\mathcal{O}}(x_k;  y_{n_{1}+ \cdots
n_{k-1} +1},\cdots , y_{n_{1}+ \cdots n_{k-1} + n_{k}}))
\end{array}$ \\
\item There is an identity element $ 1_{1}:= 1_{\mathcal{O}}\in \mathcal{O}(1)$ such that
$$\gamma^{\mathcal{O}}(x; \underbrace{1_{\mathcal{O}},\cdots ,
1_{\mathcal{O}}}_{k \hspace{2mm} \mbox{times}})= x
=\gamma^{\mathcal{O}}(1_{\mathcal{O}}; x), \quad \mbox{for all}
\quad x\in \mathcal{O}(k).$$
\item $\Sigma$-equivariant(\cite{May}).
\end{enumerate}}
\end{De}
\begin{De}
{\it Let $\mathcal{O}$ and $\mathcal{O'}$ be two $\Sigma$-operads
with respective composition products $\gamma^{\mathcal{O}}$ and
$\gamma^{\mathcal{O}^{\prime}}$ and respective associated unities
$1_{\mathcal{O}} \hspace{2mm}\mbox{and}\hspace{2mm}
1_{\mathcal{O'}}.$ A morphism of operads
$\mathcal{O}\stackrel{f}\longrightarrow\mathcal{O'}$ is a collection
$\{f_n:\mathcal{O}(n)\longrightarrow\mathcal{O'}(n)\}_{n\geq 0}$ of
$\sum$-$\mathbb{K}$-vector space homomorphisms such that:
\begin{enumerate}
\item $f(1_{\mathcal{O}})= 1_{\mathcal{O'}}$;
\item $ f_j(\gamma^{\mathcal{O}}(x_0\otimes x_1 \otimes...\otimes x_n))= \gamma^{\mathcal{O}^{\prime}}(f_n(x_0)\otimes f_{i_1}(x_1) \otimes...\otimes f_{i_n}(x_n))$ with $j=i_1 +i_2 +...+i_n;$
\item $f_n(x\ast\sigma)= f_n(x)\ast\sigma,$ with $\gamma(x\ast\sigma; a_1, \cdots ,a_n)= \gamma (x; a_{\sigma^{-1}(1)}, \cdots ,a_{\sigma^{-1}(n)}),$
for some $x \otimes a_{k_1} \otimes \cdots \otimes a_{k_n}\in
\mathcal{O}(n) \otimes \mathcal{O}(k_{1}) \otimes \cdots \otimes
\mathcal{O}(k_{n})$ and $\sigma \in \Sigma_n.$
\end{enumerate}}
\end{De}
\indent The following observations will be useful in the sequel.
\begin{Rem}
\begin{enumerate}
\item Equivalently a $\Sigma$-operad can also be defined by the so called partial
compositions:
$$\begin{array}{ccc}
 \mathcal{O}(m) \otimes \mathcal{O}(n)& \stackrel{\circ_{i}}\longrightarrow & \mathcal{O}(m + n - 1)\qquad m\geq i\geq 1\\
x \otimes y  & \longmapsto & x \circ_{i} y\\
\end{array}$$
satisfying some properties(see [\cite{Apurba}, \cite{L. P.}] for explicit axioms)\\
The two definitions are related as follows:
\begin{enumerate}
\item[(1-i)] $x \circ_{i} y = \gamma^{\mathcal{O}}(x; \overbrace{id, \cdots,\underbrace{y}_{i} \cdots  id}^{m-tuple}, );  \qquad  m\geq i\geq 1.$
\item[(1-ii)] $\gamma^{\mathcal{O}}(x;  y_{1},\cdots , y_{m})= (\cdots(((x \circ_{m} y_{m}) \circ_{m-1} y_{m-1})\cdots ) \circ_{1} y_{1}$
\end{enumerate}
\item An operad, $\mathcal{O},$ is said to be \textbf{multiplicative} if there exists a two-inputs element $m\in\mathcal{O}(2)$
such that $m\circ_1m=m\circ_2m$.
\item An operad $\mathcal{O}$ is said to be \textbf{unitary} if the unit morphism $\eta:\mathbb{K}\longrightarrow\mathcal{O}(1)$
is an isomorphism in which case the unit element of $\mathcal{O}$
will be denoted: $1_{1}= 1_{\mathcal{O}}= \eta(1_{\mathbb{K}}) \in
\mathcal{O}(1).$
\item An operad $\mathcal{O}$ is said to be \textbf{connected} if $\mathcal{O}(0)$  is isomorphic to the ground field $\mathbb{K}$(see \cite{L. P.}, 1.2.2).
To avoid confusion in the sequel, we will denote by $1_{0}$ the unit
in $\mathcal{O}(0)$
\item Let $\mathcal{O}$ be a connected operad. It has been shown in
\cite{L. P.} that for all $S\subset [n]= \{1, 2, \cdots, n \},$ with
cardinality $l<n$, $\mathcal{O}$ is equipped with a degeneracy map
$\mid_S$ defined as follows:
\begin{center}
$\begin{array}{ccc}
  \mid_S:\mathcal{O}(n) & \longrightarrow & \mathcal{O}(l) \\ p& \mapsto
& p\mid_S=p(x_1,...,x_n),
\end{array},$
where $ x_i= \left\{\begin{array}{lll}
   1_{1} \quad \mbox{if} \quad i\in S \\
   1_{0}  \quad \mbox{if not}
  \end{array}
\right.$
\end{center}
\item Any connected operad is \textbf{augmented}, that is: there exists a $\mathbb{K}$-linear map $\epsilon_j:\mathcal{O}(j)\longrightarrow\mathbb{K}.$
\item The collection of connected operads forms a sub category of the category of operads whose homomorphisms are homomorphisms of
operads $\mathcal{O}\stackrel{\psi}\longrightarrow \mathcal{P}$
satisfying: $\psi(1_{0}^{\mathcal{O}})= 1_{0}^{\mathcal{P}},$ with
$1_{0}^{\mathcal{O}} \in \mathcal{O}(0)$ and $1_{0}^{\mathcal{P}}\in
\mathcal{P}(0).$
\end{enumerate}
\end{Rem}
\indent The fundamental example is the endomorphism operad denoted
here by $ \mathcal{L}_{A}:= End_{A},$ for an object $A$ in the
category of $\mathbb{K}$-vector spaces.
\begin{Ex}  For a given associative unitary $\mathbb{K}$-algebra $A$ with multiplication
$\mu_{A}$ and unit $\eta_{A},$ the operad $\mathcal{L}_{A}$ of
multilinear homomorphisms also called endomorphism operad(see
\cite{11}, 1.2.9) defined by: for all $n\geq 0,$
$\mathcal{L}_{A}(n):= Hom_{\mathbb{K}}(A^{\otimes n}, A)$, is a non
connected unitary multiplicative $\Sigma$-operad since
$\mathcal{L}_{A}(0)\cong A\neq \mathbb{K}$. However one can consider
the sub operad $\mathcal{L}^{\prime}_{A}:=
(\mathcal{L}^{\prime}_{A}(n))_{n\geq 0}$ of the operad
$\mathcal{L}_{A}$ such that
$\mathcal{L}^{\prime}_{A}(n)=\left\{\begin{array}{ll}
    \mathcal{L}^{\prime}_{A}(n)\quad \mbox{for all} \quad n\geq 1, \\
   \mathbb{K}\qquad \mbox{if} \quad n= 0
\end{array}
\right.$ \\
Thus the operad  $\mathcal{L}^{\prime}_{A}$ is a connected
multiplicative $\Sigma$-operad with
\begin{enumerate}
 \item[$\bullet$] The associated
operadic composition $\gamma^{\mathcal{L}^{\prime}_{A}},$ the
substitution of the values of $n$ operations in a $n$-ary operation
as inputs.
\item[$\bullet$] Its associated multiplication, $m:= \mu_{A}\in \mathcal{L}^{\prime}_{A}(2)$
\item[$\bullet$] The connectedness of the operad $\mathcal{L}^{\prime}_{A}$ comes from
its definition.
\item[$\bullet$] Its associated  unit element is the identity map $A\stackrel{id_{A}}\longrightarrow A.$
\end{enumerate}
\end{Ex}

\noindent{\bf Notations}\\
 Let $\mathcal{O}= (\mathcal{O})_{n\geq 0}$ be an operad and denote by:\\
$\bullet$ $\mathcal{O}[1]= (\mathcal{O}[1](n)= \mathcal{O}(n-1))_{n\geq 1},$ the suspension of the operad $\mathcal{O};$\\
$\bullet$ $\mbox{deg}(x),$ the degree of an element $x\in
\mathcal{O}$
 so that $\mbox{deg}(x)= n$ if $x\in
\mathcal{O}(n),$ for a given $n\in \mathbb{N};$\\
$\bullet$ $\mid x \mid= \mbox{deg(x)}-1,$ the degree of the
suspension of $x\in \mathcal{O}(n).$

\section{Simplicial Structures on Operad category}
Unless otherwise stated, we consider in this section the monoidal
symmetric category of $\mathbb{K}$-vector spaces,
$\mathcal{E}_{\mathbb{K}-vect}$, and make the following assumptions:
$m \circ_{1} 1_{0} = 1_{\mathcal{O}} = m\circ_{2} 1_{0},$ for a
given  connected multiplicative operad $\mathcal{O}$ with
multiplication $m$ and  $1_{0}\in \mathcal{O}(0).$
\subsection{Construction of simplicial operads}
Let $\mathcal{O}$ be a  unitary connected multiplicative operad
with multiplication $m$ and  $1_{0}\in \mathcal{O}(0).$ Define the
following $\mathbb{K}$-vector space maps:
\begin{enumerate}
\item For all $n\in \mathbb{N}-\{0\} $ and $ 1\leq i\leq n,$ set
$$\begin{array}{ccc}
D_{i}^{n}: \mathcal{O}(n)& \longrightarrow & \mathcal{O}(n+1)
\\ p & \mapsto & D_{i}^{n}= p\circ_{i} m
\end{array}$$
 and $D^0=\eta:\mathcal{O}(0)\cong\mathbb{K}\longrightarrow\mathcal{O}(1) $
such that $ D^0(1_{0})=\eta(1_{0})=1_{1}\in\mathcal{O}(1).$ \\

\item  Once more for all
$n\in \mathbb{N} - \{0 \}$, $1\leq i\leq n$,\\ set $S_i= [n]-\{ i
\}= \{1, 2, \cdots, i-1, \hat{i}, i+1, \cdots, n\}$, where $\hat{i}$
means that the natural number $i$ has been omitted and for $n> 1,$
define
$$\begin{array}{ccc} F_i^n= \mid_{S_i}: \mathcal{O}(n)&
\longrightarrow & \mathcal{O}(n-1)\\
 p & \mapsto & \begin{array}{cc}
F_i^n(p)=p\mid_{S_i} &= p\circ_i 1_{0}
                    \end{array}
\end{array}$$

 and $F^1=\mid_{\emptyset}:\mathcal{O}(1)\longrightarrow\mathcal{O}(0)$
such that for $1_{\mathcal{O}}\in\mathcal{O}(1)$,
$F^1(1_{\mathcal{O}})=1_{\mathcal{O}}\mid_{\emptyset}=1_{0}\in\mathcal{O}(0).$
\end{enumerate}

\indent  One can verify by straightforward calculations that the
$\mathbb{K}$-linear maps $D_i^n$ and $F_i^n$ respectively called
degeneracy map and face map are subject to the following properties:
\begin{enumerate}
 \item[(i)]  $$F_iF_j=F_{j-1}F_i \quad {\mbox{if}}\quad i<j \qquad(2.3)$$
 \item[(ii)] $$D_iD_j=D_{j+1}D_i \quad {\mbox{if} }\quad i\leq j \qquad(2.4)$$
\item[(iii)] $$
 F_iD_j =\left\{\begin{array}{lll}
    D_{j-1}F_i\quad \mbox{if} \quad i<j \\
  id  \qquad \mbox{if} \quad i\in\{j,j+1\} \qquad (2.5)\\
  D_jF_{i-1}  \quad \mbox{if} \quad i>j+1
\end{array}
\right.$$
\item[$(iv)$]$$F_iF_j=F_jF_{i+1} \quad \mbox{if} \quad i\geq j \quad (2.3')$$
\item[$(v)$]$$D_iD_j=D_jD_{i-1} \quad \mbox{if} \quad j<i \quad \quad (2.4')$$
\end{enumerate}

\begin{De}{\it A simplicial operad is a unitary connected multiplicative
operad endowed with face homomorphisms, $F_i,$ and degeneracy
homomorphisms, $D_i.$}
\end{De}

\begin{Rem}
 {\it Since the above construction of simplicial structure is functorial then  for any homomorphism
  of simplicial operads $\lambda:\mathcal{O} \rightarrow \mathcal{O'},$ it is obvious that $F_i\lambda_n=\lambda_{n-1}F_i,$ $n\geq 1$
and $D_i\lambda_n=\lambda_{n+1}D_i,$  $n\geq 0$.}\\
\end{Rem}

\subsection{A Differential graded coalgebra structure on simplicial operads}
\indent Now let $\mathcal{O}$ be a simplicial operad. The $-1$
degree $\mathbb{K}$-linear map
$\mathcal{O}\stackrel{\partial}\longrightarrow \mathcal{O}$  defined
for all $n\geq 1,$ as follows:
$$\begin{array}{ccccc}\partial_n:\mathcal{O}(n)&  \longrightarrow &\mathcal{O}(n-1)& &\\
  & & & & (2.6)\\ x_n&  \mapsto & \partial_{n}(x_{n})= &
  \sum\limits_{i=1}^{n}(-1)^iF_i^nx_n=\sum\limits_{i=1}^{n}(-1)^ix_n\circ_i1_0
  &
 \end{array}$$
 leads us to the result below:
 \begin{Pro} {\it If $\mathcal{O}$ is a
simplicial operad, Then the family of pairs
$(\mathcal{O}(n),\partial_n)_{n\in\mathbb{N}}$ is a
$\mathbb{K}$-vector space chain complex with boundary operator
$\partial$.}
\end{Pro}
 \begin{proof}\\ Let
$\mathcal{O}$ be a simplicial operad. For $n\in\mathbb{N}$
 and $x_{n+1}\in\mathcal{O}(n+1)$. One has:
 $$\partial_n\partial_{n+1}x_{n+1}=\sum\limits_{i=1}^{n}(-1)^i\sum\limits_{j=1}^{n+1}(-1)^jF_i^nF_j^{n+1}x_{n+1}=
 \sum\limits_{i=1}^{n}(-1)^i\sum\limits_{j=1}^{n+1}(-1)^j(x_{n+1}\circ_j1_0)\circ_i1_0.$$
  Considering the disjunction cases on the variables $i\in \{1,
\cdots, n\}$ and using the relation (2.3'), we successively obtain:
\begin{enumerate}
\item[$\bullet$] for $i=1:   F_1^nF_2^{n+1}-F_1^nF_2^{n+1}+F_1^nF_3^{n+1}+...+(-1)^{n+2}F_1^nF_{n+1}^{n+1}$,
\item[$\bullet$] for $i=2:  -F_1^nF_3^{n+1}+F_2^nF_3^{n+1}-F_2^nF_3^{n+1}+F_2^nF_4^{n+1}+F_2^nF_5^{n+1}+...+(-1)^{n+1}F_2^nF_{n+1}^{n+1}$,
\item[$\bullet$] for $i=3:
F_1^nF_4^{n+1}-F_2^nF_4^{n+1}+F_3^nF_4^{n+1}-F_3^nF_4^{n+1}+F_3^nF_5^{n+1}+...+(-1)^{n+2}F_3^nF_{n+1}^{n+1}$.\\
 $\cdot$\\
 $\cdot$\\
 $\cdot$
\item[$\bullet$] for $i=n-1:  (-1)^{n-1}(-F_1^nF_n^{n+1}+F_2^nF_n^{n+1}+...+(-1)^{n-1}F_{n-1}^nF_n^{n+1}+(-1)^nF_{n-1}^nF_n^{n+1}+(-1)^{n+1}F_{n-1}^nF_{n+1}^{n+1})$
\item[$\bullet$]  for $i=n:  (-1)^{n}(-F_1^nF_{n+1}^{n+1}+F_2^nF_{n+2}^{n+1}+...+(-1)^{n}F_{n}^nF_{n+1}^{n+1}+(-1)^{n+1}F_{n}^nF_{n+1}^{n+1}$ .
\end{enumerate}
Taking the sum of all these expressions, we obtain the following
identity:
$$\partial_n\partial_{n+1}x_{n+1}=0, \quad \mbox{for all}\quad n\in \mathbb{N}.$$\\ Notice that this result can also
be deduced from the property (2.3) of the above defined simplicial
maps.
\end{proof}\\\\
 \indent In the sequel we are going to define a coalgebra structure on a
simplicial operad through the Alexander-Whitney map.\\
 \indent Let $\mathcal{O}$ be a simplicial operad. For a given non zero
integer $n$, consider the maps $\widetilde{F}^{n-j}_{n}:
\mathcal{O}(n)\longrightarrow \mathcal{O}(j)$ and
$\widetilde{F}_1^{j}: \mathcal{O}(n)\longrightarrow
\mathcal{O}(n-j)$ such that for all $x\in \mathcal{O}(n),$ one has
respectively:
$$\widetilde{F}_{n}^{n-j}x = \left\{\begin{array}{ll}F_{j+1}^{j+1}\circ
 F_{j+2}^{j+2}\circ \cdots \circ F_{n-1}^{n-1}\circ F_{n}^{n}x, \quad if \quad 0 \leq j<n,\\
x,  \quad if \quad j= n \end{array} \right.$$
$$\widetilde{F}_1^{j}x =
\left\{\begin{array}{ll}F_{1}^{n-j+1}\circ
 F_{1}^{n-j+2}\circ \cdots \circ F_{1}^{n-1}\circ F_{1}^{n}x, \quad if \quad 0< j\leq n,\\
x,  \quad if \quad j= 0 \end{array} \right.$$.\\
From the maps constructed above follow the graded homomorphisms:
\begin{enumerate}
\item For all $n\geq 0,$ $\begin{array}{ccc}
\Delta_{\mathcal{O}}:\mathcal{O}(n)& \longrightarrow &
(\mathcal{O}\otimes\mathcal{O})(n)=
\bigoplus\limits_{p+q=n}\mathcal{O}(p)\otimes\mathcal{O}(q)\\ x &
\longmapsto &
\sum_{j=0}^{n}{\widetilde{F}_{n}^{n-j}x\otimes\widetilde{F}_1^jx}
\end{array}$
\item $\epsilon:\mathcal{O}\longrightarrow\mathbb{K}$ is a zero homomorphism in non zero degree and since the operad  $\mathcal{O}$ is connected,
$\epsilon$ is the identity homomorphism in degree zero.
\end{enumerate}
These homomorphisms lead us the result below:
\begin{Pro}
Let $\mathcal{O}$ be any simplicial operad  in the category of
vector spaces, the triple $((\mathcal{O}(n),\partial),
\Delta_{\mathcal{O}},\epsilon)$ is a counitary differential graded
coalgebra.
\end{Pro}
\begin{proof}
 To show that $(\mathcal{O}(n), \Delta_{\mathcal{O}})$ a graded coalgebra it is enough to prove that the coproduct $\Delta_{\mathcal{O}}$ is coassociative
 and this is equivalent to verify that the following diagram  commutates\\
$$\xymatrix{
 \mathcal{O}\ar[d]_{\Delta_{\mathcal{O}}} \ar[r]^{\Delta_{\mathcal{O}}}& \mathcal{O}\otimes\mathcal{O} \ar[d]^{id\otimes\Delta_{\mathcal{O}}}  &\\
\mathcal{O}\otimes\mathcal{O} \ar[r]^{\Delta_{\mathcal{O}}\otimes id}& (\mathcal{O}\otimes\mathcal{O})\otimes\mathcal{O} &\\
}.$$ That is: $(id\otimes\Delta_{\mathcal{O}})\circ\Delta_{\mathcal{O}}=(\Delta_{\mathcal{O}}\otimes id)\circ\Delta_{\mathcal{O}}$.
In fact one has for $x\in\mathcal{O}(n)$\\
\begin{align*}
    [(id\otimes\Delta_{\mathcal{O}})\circ \Delta_{\mathcal{O}}](x)&=(id\otimes\Delta_{\mathcal{O}})
    (\sum\limits_{i=0}^n\widetilde{F}_{n}^{n-i}x\otimes\widetilde{F}_1^{i}x)\\
    &=\sum\limits_{i=0}^n\widetilde{F}_{n}^{n-i}x\otimes\Delta(\widetilde{F}_1^{i}x)\\
    &=\sum\limits_{i=0}^n\widetilde{F}_{n}^{n-i}x\otimes\sum\limits_{j=0}^{n-i}
    \widetilde{F}_{n}^{n-i-j}\widetilde{F}_1^{i}x\otimes\widetilde{F}_1^{j}
    \widetilde{F}_1^{i}x\\
    &=\sum\limits_{i=0}^n\sum\limits_{j=i}^{n}\widetilde{F}_{n}^{n-i}x\otimes
    \widetilde{F}_{n}^{n-j}\widetilde{F}_1^{i}x\otimes\widetilde{F}_1^{j-i}
    \widetilde{F}_1^{i}x\\
    &=\sum\limits_{i=0}^n\sum\limits_{j=i}^{n}A_{i,j}.
\end{align*}

Similarly we have\\
\begin{align*}
    [(\Delta_{\mathcal{O}}\otimes id)\circ \Delta_{\mathcal{O}}](x)&=(\Delta_{\mathcal{O}}\otimes id)
    (\sum\limits_{i=0}^n\widetilde{F}^{n-i}x\otimes\widetilde{F}_1^{i}x)\\
    &=\sum\limits_{i=0}^n\Delta_{\mathcal{O}}(\widetilde{F}^{n-i}x)\otimes\widetilde{F}_1^{i}x\\
    &=\sum\limits_{i=0}^n\sum\limits_{j=0}^{i}\widetilde{F}^{i-j}\widetilde{F}^{n-i}x\otimes
    \widetilde{F}_1^{j}\widetilde{F}^{n-i}x\otimes\widetilde{F}_1^{i}x\\
    &=\sum\limits_{i=0}^n\sum\limits_{j=0}^{i}B_{i,j}.
\end{align*}
The sums above have exactly
$\sum\limits_{i=0}^{n}n-i+1=\sum\limits_{i=1}^{n+1}i=\frac{(n+1)(n+2)}{2}$
terms and  these terms satisfy the following equalities.
\begin{enumerate}
    \item $A_{i,i}=B_{i,i},\quad 0\le i\le n$.
    \item $A_{0,j}=B_{n,j}, \quad 1\le j\le n-1$.
    \item $A_{i,j}=B_{j,i}, \quad i<j$.
\end{enumerate}
Thus the coproduct $\Delta_{\mathcal{O}}$ is coassociative.\\
 Furthermore one can verify by straightforward computation that $\partial$ is a
coderivation with respect to the coproduct $\Delta_{\mathcal{O}}.$\\
 So a simplicial operad has a co-unitary differential graded
coalgebra structure.
\end{proof}

\subsection{Link between the composition product and the
face(respectively degeneracy) map}
 Let $s\in \mathbb{N},$
$(t_{r})_{r=1}^{s}$ and $(j_{r})_{r=1}^{s}$ two s-tuples whose
components are in $\mathbb{N}$ such that $(1\leq j_{r}\leq
t_{r})_{r=1}^{s}.$
\begin{enumerate}
\item From the face homomorphism define the following
$\mathbb{K}$- linear homomorphisms:\\
\begin{enumerate}
  \item[1.a-]  \begin{center} $\begin{array}{ccc} F_{(j_{1},\cdots, j_{s})}^{(t_{1}, \cdots,
t_{s})}
:\mathcal{O}(t_{1})\otimes\cdots \otimes\mathcal{O}(t_{s})&\longrightarrow & \mathcal{O}(t_{1}-1)\otimes\cdots \otimes\mathcal{O}(t_{s}-1)\\
p_{1}\otimes\cdots \otimes p_{s}& \mapsto &
(-1)^{\varepsilon_{s}}F_{j_{1}}^{t_{1}}(p_{1})\otimes\cdots\otimes
F_{j_{s}}^{t_{s}}(p_{s}),
\end{array}$
 $\mbox{where}\quad \varepsilon_{s}= \sum\limits_{r=1}^{s}
 t_{r}(s-r).$
 \end{center}

  \item[1.b-] Set $t_0= 0,$
 $$ \begin{array}{ccc} F_{j_{1} + j_{2}\cdots + j_{s}}^{t_{1} +
t_{2}
 \cdots + t_{s}}:\mathcal{O}(t_1+...+t_s)&\longrightarrow &\mathcal{O}(t_{1} +
t_{2} + \cdots + t_{s}-s)\\  p &\longmapsto & (((\cdots((\cdots
(p\circ_{(j_{s}+ \sum\limits_{r=0}^{s-1}t_{r})}1_{0}
)\cdots)\circ_{(j_{l}+
\sum\limits_{r=0}^{l-1}t_{r})}1_{0})\cdots)\\
& &  \circ_{t_{1}+j_{2}}1_{0}) \circ_{j_{1}}1_{0})
\end{array}$$
\end{enumerate}
Notice that although the two maps $F_{(j_{1},\cdots,
j_{s})}^{(t_{1}, \cdots, t_{s})}$ and $F_{j_{1} + j_{2}\cdots +
j_{s}}^{t_{1} + t_{2}
 \cdots + t_{s}}$ defined above are homomorphisms of degree $-s$,
 the first one
 acts on $s$ operations $(p_{i})_{i=1}^{i=s}$ with respective
 inputs $(t_i)_{i=1}^{i=s}$ while the second acts on only one operation with total
 inputs
 $\sum\limits_{i=1}^{s}t_i$

\item On the other hand, one obtains from the degeneracy maps the
$\mathbb{K}$-linear maps:\\

 \begin{enumerate}
  \item[2.a-] \begin{center}
$D_{(j_{1},\cdots, j_{s})}^{(t_{1}, \cdots, t_{s})}
:\mathcal{O}(t_{1})\otimes\cdots \otimes\mathcal{O}(t_{s})\longrightarrow\mathcal{O}(t_{1}+1)\otimes\cdots \otimes\mathcal{O}(t_{s}+1)$\\
\hspace{20cm} $p_{1}\otimes\cdots \otimes p_{s}\mapsto
(-1)^{\varepsilon_{s}}D_{j_{1}}^{t_{1}}(p_{1})\otimes\cdots
\otimes D_{j_{s}}^{t_{s}}(p_{s}),$\\
 $\mbox{where}\quad
\varepsilon_{s}= \sum\limits_{r=1}^{s} t_{r}(s-r).$
\end{center}

\item[2.b-] Set once more $t_0= 0,$
\begin{center} $ \begin{array}{ccc} D_{j_{1} + j_{2}\cdots
+ j_{s}}^{t_{1} + t_{2} \cdots + t_{s}}:\mathcal{O}(t_1+...+t_s) &
\longrightarrow & \mathcal{O}(t_{1} +
t_{2} + \cdots + t_s+ s)\\
 p&\longmapsto &
(((\cdots((\cdots (p\circ_{(j_{s}+ \sum\limits_{r=0}^{s-1}t_{r})}m
)\cdots)\circ_{(j_{l}+
\sum\limits_{r=0}^{l-1}t_{r})}m)\cdots)\\
& &  \circ_{t_{1}+j_{2}}m) \circ_{j_{1}}m).
\end{array}$
\end{center}

Similarly $D_{(j_{1},\cdots, j_{s})}^{(t_{1}, \cdots, t_{s})}$ and
$D_{j_{1} + j_{2}\cdots + j_{s}}^{t_{1} + t_{2}
 \cdots + t_{s}}$ defined above are two homomorphisms of degree $s$
 acting respectively on $s$ operations $(p_{i})_{i=1}^{i=s}$ with respective
 inputs $(t_i)_{i=1}^{i=s}$ and one operation with total inputs
 $\sum\limits_{i=1}^{s}t_i.$

\end{enumerate}
\end{enumerate}

Thus the composition product is respectively related to the face map
and the degeneracy map through the following identities:\\
 $$F_{j_{1}+ j_{2}+ \cdots +j_{s}}^{t_{1} + t_{2} + \cdots +
t_{s}}\circ\gamma= (-1)^{\varepsilon_{s}}\gamma\circ(Id\otimes
F_{(j_{1},\cdots,j_{s})}^{(t_{1},\cdots,t_{n})}),$$

$$D_{j_{1}+ j_{2}+ \cdots +j_{s}}^{t_{1} + t_{2} + \cdots +
t_{s}}\circ\gamma= (-1)^{\varepsilon_{s}}\gamma\circ(Id\otimes
D_{(j_{1},\cdots,j_{s})}^{(t_{1},\cdots,t_{n})}).$$

\section{Simplicial differential graded right brace algebra structure.}
In this section we are going to define on a simplicial operad a new
brace algebra structure from which we will obtain an analogous of
homotopy G-algebra structure.\\\\
 \indent{\bf 4.1-} Consider $\mathcal{O} = \bigoplus _{k\geq 0}
\mathcal{O}(k),$ a simplicial operad in the category of the
$\mathbb{K}$-vector spaces, $\mathcal{E}_{\mathbb{K}-vect},$ and for
all $s\in\mathbb{N},$ set
$$(\mathcal{O}^{\otimes n})(s)= \bigoplus _{s_1+s_2+\cdots + s_n= s}
\mathcal{O}(s_1)\otimes \mathcal{O}(s_2)\otimes\cdots \otimes
\mathcal{O}(s_n).$$ If $p \in \mathcal{O}(r)$, we recall that
$deg(p) = r$  and $|p| = r-1$ denote respectively the degree of $p$ and the degree of its suspension. \\
A right brace structure on $\mathcal{O}$ is a collection of
multilinear operations defined by:
\begin{center} $\mathcal{O}\otimes\mathcal{O}^{\otimes n}\longrightarrow\mathcal{O}$, $n\geq 1$\\
 $p\otimes(q_1\otimes q_2\otimes\cdots\otimes q_n)\mapsto
 p\{q_1,q_2,\cdots ,q_n\} $
\end{center} such that for $p, q_1,q_2,\cdots,q_n \in\mathcal{O},$ $$p\{q_1,q_2,\cdots,q_n\}:=\sum\limits_{}^{}(-1)^{\epsilon}\gamma(p; 1_{\mathcal{O}},...,
1_{\mathcal{O}},q_1,1_{\mathcal{O}},...,1_{\mathcal{O}}, q_2,
1_{\mathcal{O}},...,1_{\mathcal{O}},
q_n,1_{\mathcal{O}},...,1_{\mathcal{O}}),$$
 where the summation runs
over all possible substitutions of $q_1, q_2,\cdots, q_n$ into $p$
in the prescribed order and $\epsilon:=\sum\limits_{j=1}^{n}\mid
q_j\mid(deg(p\{q_1,\cdots,q_n\})-i_j))$ , $i_j$ being the total
number of inputs in front of $q_j$.
 Thus the right braces $p\{q_1,q_2,\cdots,q_n\}$ are homogeneous of degree $-n,$ i.e.
 $deg(p\{q_1,q_2,\cdots,q_n\})= deg(p) + \sum\limits_{j=1}^{n}\mid q_j\mid $ with
 the following conventions: $p\circ q := p\{q\}$ and $p\{\}:= p,$ for
 all $p, q \in \mathcal{O}.$\\
\begin{De}
An operad $\mathcal{O}$ is said to be endowed with a simplicial
right brace algebra structure if it is a simplicial operad together
with a right brace algebra structure.
\end{De}
\indent{\bf 4.2-} As in the case of Gerstenhaber-Vornov's brace
structure,
 observe that the right brace structure on an operad $\mathcal{O}$ induces  an
associative product called dot-product and defined as follows:
$$p\bullet q =(-1)^{(\mid q\mid +1)(\mid p\mid +1)}\gamma^{\mathcal{O}}(m;p,q), \quad \mbox{for}\quad  p, q\in
\mathcal{O}.$$ But this product satisfies neither the distributivity
and nor the higher homotopy properties of Gerstenhaber-Voronov
[\cite{VG}, \cite{G. V.}]. Thus, to find  an analogous of their
homotopy G-algebra(see [\cite{VG}, \cite{G. V.}] for explicit
definition), we introduce on an operad $\mathcal{O}$ equipped with a
right brace structure, a new associative product called odot-product
by setting:
 $$ p\odot q := (-1)^{s(r-1)}m\{p,q\} \quad \mbox{for} \quad p\in\mathcal{O}(r)\quad \mbox{and} \quad q\in\mathcal{O}(s).$$
 Notice that this odot-product satisfies among other things, the distributivity
and the higher homotopies properties.
 Moreover the two products are related by: $$p\bullet q = (-1)^{rs}p\odot q.$$
 \indent Consider a simplicial operad $\mathcal{O}$ endowed a right brace structure.  There exists a
 coboundary operator defined by: for all $n\in \mathbb{N},$ $$\begin{array}{ccc} \mathcal{O}(n) & \stackrel{d}\longrightarrow &
\mathcal{O}(n+1)\\ p& \longmapsto & d(p)= m\circ
p+\sum\limits_{i=1}^{\mid p\mid +1}(-1)^iD_ip \end{array}.$$ This
differential operator is not a derivation with respect to the
dot-product, in contrast we have the statement below:
\begin{Pro}
If $\mathcal{O}= \bigoplus_{n\geq 0}\mathcal{O}(n)$ a simplicial
right brace algebra then the following assertions hold:
\begin{enumerate}
 \item  $((\mathcal{O},d), \odot)$ is a differential graded algebra.
 \item  $(\mathcal{O}, \partial, d)$ is a bicomplex with total differential $D= d +
 \partial.$
 \end{enumerate}
\end{Pro}
\begin{proof}
\begin{enumerate}
\item
To show that $((\mathcal{O},d), \odot)$ is a differential graded
algebra, it is enough to establish that $d$ is a derivation with
respect to the odot-product, $\odot.$\\
So let $p\in\mathcal{O}(r)$ and $q\in\mathcal{O}(s)$.\\ By using the
identities $d(p\bullet q)= (-1)^{s}d(p)\bullet q + p\bullet d(q),$
$p\odot q=(-1)^{rs}p\bullet q$, we have
\begin{align*}
d(p\odot q)&=(-1)^{rs}d(p\bullet q)\\
    &=(-1)^{rs}[(-1)^{s}d(p)\bullet q + p\bullet d(q)]\\
    &=(-1)^{rs}[(-1)^s(-1)^{s(r+1)}d(p)\odot q +(-1)^{r(s+1)}p\odot d(q)]\\
    &=d(p)\odot q+(-1)^r p\odot d(q).
\end{align*}
So d is a derivation with respect to the odot-product $\odot$.

\item  Hereafter we are going to show that $\mathcal{O}$ is a
bicomplex with total differential $D= d + \partial.$\\ Since $d$ and
$\partial$ are differential operators, it is enough to check that
the two differential operators commute in graded sense. So let
$\mathcal{O}$ be simplicial operad with multiplication $m,$ units
$1_{\mathcal{O}}\in \mathcal{O}(1)$ and $1_{0}\in \mathcal{O}(0).$
We recall the following properties for further use. For all
$f\in\mathcal{O}(t)$, $g\in\mathcal{O}(n)$,
 $h\in\mathcal{O}(l)$, $m\in\mathcal{O}(2),$ one has:\\
$\left\{\begin{array}{lll}
(f\circ_i g)\circ_{j+i-1}h=f\circ_i(g\circ_j h);\qquad (a)\\
(f\circ_i g)\circ_{j+n-1}h=(f\circ_j h)\circ_i g, i<j;\qquad (b)\\
m\circ_2 1=e_{\mathcal{O}}=m\circ_1 1;\qquad (c)\\
 F_iD_j =\left\{\begin{array}{lll}
    D_{j-1}F_i\quad \mbox{if} \quad i<j \\
  id  \qquad \mbox{if} \quad i\in\{j,j+1\} \qquad (d)\\
  D_jF_{i-1}  \quad \mbox{if} \quad i>j+1
\end{array}
\right.
\end{array}
\right.$ \\ Moreover for all $f\in\mathcal{O}(n),$
$d(f)=m\{f\}+(-1)^n f\{m\} $ and
 $\partial(f)= \sum\limits_{i=1}^{n}(-1)^i f\circ_i 1_{0}.$ So we
 have:
$$\begin{array}{lll}
(\partial\circ d)(f)&=\partial(m\{f\}+(-1)^nf\{m\});\\
                   &=\sum\limits_{i=1}^{n+1}(-1)^i(m\{f\})\circ_i 1_{0}+(-1)^{n}\sum\limits_{i=1}^{n+1}(-1)^i(f\{m\})\circ_i 1_{0}\\
                   &=(-1)^{n-1}\sum\limits_{i=1}^{n+1}(-1)^i(m\circ_1f)\circ_i 1_{0} +\sum\limits_{i=1}^{n+1}(-1)^i(m\circ_2f)\circ_i 1_{0}\\
                   &+\sum\limits_{i=1}^{n+1}\sum\limits_{j=1}^{n}(-1)^{i+j}(f\circ_jm)\circ_i 1_{0}\\
                   &=(-1)^{n-1}\sum\limits_{i=1}^{n}(-1)^i(m\circ_1f)\circ_i 1_{0} + (m\circ_1f)\circ_{n+1}1_{0} -(m\circ_2 f)\circ_1 1_{0}\\
                  &+\sum\limits_{i=2}^{n+1}(-1)^i(m\circ_2f)\circ_i 1_{0} +\sum\limits_{i=1}^{n+1}\sum\limits_{j=1}^{n}(-1)^{i+j}(f\circ_jm)\circ_i 1_{0}\\
                  &=(-1)^{n-1}\sum\limits_{i=1}^{n}(-1)^i(m\circ_1 f)\circ_i 1_{0} + \underbrace{(m\circ_2 1_{0})\circ_{1}f-(m\circ_1 1_{0})\circ_1f}_{=0}\\
                  &+\sum\limits_{i=2}^{n+1}(-1)^i(m\circ_2f)\circ_i 1_{0} + \sum\limits_{i=1}^{n+1}\sum\limits_{j=1}^{n}(-1)^{i+j}(f\circ_jm)\circ_i 1_{0}\\
\end{array}$$
          $$\begin{array}{lll}
                                    &=(-1)^{n-1}\sum\limits_{i=1}^{n}(-1)^i(m\circ_1 f)\circ_i 1_{0} + \sum\limits_{i=2}^{n+1}(-1)^i(m\circ_2f)\circ_i 1_{0}\\
                                    &+(-1)^{n}\sum\limits_{i=1}^{n+1}\sum\limits_{j=1}^{n}(-1)^{i+j}(f\circ_jm)\circ_i 1_{0}\\
                                      &=(-1)^{n-1}\sum\limits_{i=1}^{n}(-1)^i\underbrace{(m\circ_1 f)\circ_i 1_{0}}_{=m\circ_1(f\circ_i 1_{0})}+\sum\limits_{i=2}^{n+1}(-1)^i\underbrace{(m\circ_2f)\circ_i 1_{0}}_{=m\circ_2(f\circ_{i-1} 1_{0})}\\
                                      &+\sum\limits_{i=1}^{n+1}\sum\limits_{j=1}^{n}(-1)^{i+j}(f\circ_jm)\circ_i 1_{0},\quad (a)\\
                                      &=(-1)^{n-1}\sum\limits_{i=1}^{n}(-1)^im\circ_1(f\circ_i 1_{0})+\sum\limits_{i=2}^{n+1}(-1)^im\circ_2(f\circ_{i-1}1_{0})\\
                                      &+\sum\limits_{i=1}^{n+1}\sum\limits_{j=1}^{n}(-1)^{i+j}(f\circ_jm)\circ_i 1_{0};\\
                                      &=(-1)^{n-1}\sum\limits_{i=1}^{n}(-1)^im\circ_1(f\circ_i 1_{0})-\sum\limits_{i=1}^{n}(-1)^im\circ_2(f\circ_{i}1_{0})\\
                                      &+\sum\limits_{i=1}^{n+1}\sum\limits_{j=1}^{n}(-1)^{i+j}(f\circ_jm)\circ_i 1_{0};\\
                                      &=(-1)^{n-1}\sum\limits_{i=1}^{n}(-1)^im\circ_1(f\circ_i 1_{0})-\sum\limits_{i=1}^{n}(-1)^im\circ_2(f\circ_{i}1_{0})\\
                                      &+\sum\limits_{i=1}^{n+1}\sum\limits_{j=1}^{n}(-1)^{i+j}F_i^{n+1}D_j^nf\\
                                      &=(-1)^{n-1}\sum\limits_{i=1}^{n}(-1)^im\circ_1(f\circ_i 1_{0})-\sum\limits_{i=1}^{n}(-1)^im\circ_2(f\circ_{i} 1_{0})\\
                                      &-\sum\limits_{i=1}^{n}\sum\limits_{j=1}^{n-1}(-1)^{i+j-1}(f\circ_i 1_{0})\circ_jm\\
                   &=-[(-1)^{n-2}\sum\limits_{i=1}^{n}(-1)^im\circ_1(f\circ_i 1_{0})+\sum\limits_{i=1}^{n}(-1)^im\circ_2(f\circ_i 1_{0})\\
                   &+\sum\limits_{i=1}^{n}\sum\limits_{j=1}^{n-1}(-1)^{i+j-1}(f\circ_i 1_{0})\circ_jm];\\
                   &=-(\partial\circ d)(f).
\end{array}$$
And since $\partial^{2}= 0= d^{2},$ then the result follows.
\end{enumerate}
\end{proof}
\\
 Hereafter, we make the following observations to highlight the differences
  and analogies between the above defined brace algebra structure and the one given by Gerstenhaber-Voronov.
 \begin{Rem}
  \begin{enumerate}
  \item The above mentioned Gerstenhaber-Voronov's brace algebra structure is
 given by: for $m\geq n\geq 1$
and $p, q_1,q_2,\cdots,q_n \in\mathcal{O},$
$$p\{q_1,q_2,\cdots,q_n\}(x_1, \cdots, x_m):=\sum\limits_{}^{}(-1)^{\epsilon}p(
x_{1},..., x_{i_1},q_1(x_{i_1 + 1}, \cdots),\cdots ,x_{i_n},
q_n(x_{i_n + 1} \cdots), \cdots, x_{m}),$$
 where $\epsilon:=\sum\limits_{j=1}^{n}\mid
q_j\mid i_j$ , that is $(-1)^{\epsilon}$ is the sign picked up by
rearranging the sequence of operations on the left-hand side of the
formula into the sequence of terms on the right-hand side, in
accordance to the usual sign convention.\\
  Observe that unlike the Gerstenhaber-Voronov's definition in which they consider the number of total inputs in front of a given $q_j$ from the left to the
 right to fix the sign of the brace, in our case this sign is defined by using the number of total inputs in front of
  the same $q_j$ from the right to the left except the inputs of $q_j,$ hence the terminology right
  brace.\\

\item The right brace also satisfies the associativity property known as higher pre-Jacobi identity and recalled below:
 $$ \begin{array}{cc} x\{ x_1, x_2, \cdots, x_{s} \}\{ y_1, y_2, \cdots, y_{r} \}=& \sum\limits_{r\geq i_{r}\geq i_{r-1}\geq
 \cdots \geq i_{1}\geq 0}(-1)^{\varepsilon}x\{y_1, \cdots y_{i_{1}},
 x_{1}\{y_{i_{1} +1}, \cdots, y_{j_{1}}\},\\ & y_{j_{1}
 + 1}, \cdots, y_{i_{s}}, x_{s}\{y_{i_{s} +1}, \cdots, y_{j_{s}
 }\}, y_{j_{s}+ 1}  \cdots, y_{r} \}
 \end{array}$$
 where
$\varepsilon=\sum\limits_{p=1}^{s}\mid
x_p\mid(\sum\limits_{q=j_{p}+1}^{r}\mid y_q\mid)$ with $j_p$ the
position of the last $y_i$ in the right brace
$x_p\{y_{i_p+1},\cdots,y_{j_p}\}.$
\item We also specify that our dot-product as well as our odot-product are different from the Gerstenhaber-Voronov's
product, and therefore our resulting algebra is different from their
algebra.
\end{enumerate}
\end{Rem}

Hereafter we introduce on a simplicial operad equipped with a right
brace structure, an analogous of homotopy G-algebra called right
homotopy G-algebra.
\begin{De}
A right homotopy G-algebra is a differential graded algebra
$((\mathcal{O}, d), \odot )$ together  with a right brace structure
satisfying the distributive  and higher homotopy properties.

\end{De}
\begin{Th}
Let $\mathcal{O}= \bigoplus_{k\geqslant 0}\mathcal{O}(k)$ be a
simplicial graded operad together with its right brace algebra
structure. The triple $((O, \partial), \odot)$ is a differential
graded algebra.
 Moreover for every $p\in \mathcal{O}(r), (z_{j}\in
\mathcal{O}(t_{j}))_{1\leq j\leq n}^{t_{j}\geq 1},  r> n,$ the
differential operator $\partial$ satisfies the relation:
$$(\ast \ast) \quad \partial_{u}(p\{z_1,\cdots,z_n\})= \left\{ \begin{array}{ll}
 \sum\limits_{s= 1}^{n}(-1)^{\delta} p\{z_1,\cdots,(\partial z_s),\cdots,z_{n}\} \quad \mbox{if r is even} \\
  -(F_r^rp)\{z_1,\cdots,z_n\}+ \sum\limits_{s= 1}^{n}(-1)^{\delta}p\{z_1,\cdots,(\partial z_s),\cdots,z_{n}\}, \quad \mbox{else.}\\
   \end{array} \right.$$
   where $\delta=deg(p\{z_1,\cdots,z_n\})-deg(z_s)+\sum\limits_{i=1}^{s-1}\mid
   z_i\mid$ and
$u= r + (\sum\limits_{j=1}^{n}  t_{j}) -n , r> n.$
\end{Th}
\begin{proof}
To prove the first part of this result we show that the differential
$\partial$ is a derivation with respect to the
odot-product $\odot.$\\ Consider for this purpose $p\in\mathcal{O}(n),q\in\mathcal{O}(l)$, we have:\\

\begin{equation*}
\begin{aligned}
\partial_{n+l}(p\odot q)&=& \sum\limits_{i=1}^{n+l}(-1)^iF_i^{n+l}m(p,q)\\
          &=& m(\sum\limits_{i=1}^{n}(-1)^iF_i^{n}p, q)+(-1)^n m(p,\sum\limits_{i=n+1}^{n+l}(-1)^iF_i^{n+l}q)\\
          &=& m(\sum\limits_{i=1}^{n}(-1)^iF_i^{n}p, q)+(-1)^n m(p,\sum\limits_{i=1}^{l}(-1)^iF_i^{l}q)]\\
          &=& m(\partial_np,q)+(-1)^n m(p,\partial_lq)\\
          &=&(\partial_n p)\odot q+ (-1)^n p\odot\partial_lq
\end{aligned}
\end{equation*}
 Thus $\partial$ is a derivation with respect to the odot-product
$\odot$ and therefore the triple $((O, \partial), \odot)$ is a
differential graded algebra.\\
On the other hand the relation $(\ast \ast)$ is established as
follows:\\ Given $p\in \mathcal{O}(r), (z_{j}\in
\mathcal{O}(t_{j}))_{1\leq j\leq n}^{t_{j}\geq 1}.$ One has the
right brace $p\{z_{1},z_{2} \cdots z_{n}\}\in \mathcal{O}(u)$ with
$u= r + \sum\limits_{j=1}^{n} \mid z_{j} \mid , r> n$ and:
\begin{equation*}
\begin{aligned}
\partial_{u}(p\{z_1,\cdots,z_{n}\})&= \sum\limits_{i=1}^{u}(-1)^iF_i^u(p\{z_1,\cdots,z_{n}\})\\
                                &= \sum\limits_{i=1}^{u}(-1)^iF_i^u[\sum\limits_{1\leq j_1\leq \cdots \leq j_{s} \leq\cdots \leq j_{n}}(-1)^{\epsilon}p(1_{\mathcal{O}},\\
                                &\cdots,1_{\mathcal{O}},z_1,1_{\mathcal{O}},\cdots, 1_{\mathcal{O}}, z_s, 1_{\mathcal{O}},\cdots 1_{\mathcal{O}},z_{n},1_{\mathcal{O}},\cdots,1_{\mathcal{O}})]\\
                                 &= \sum\limits_{1\leq j_1\leq\cdots \leq j_{s}\leq \cdots \leq j_{n}}(-1)^{\epsilon}\sum\limits_{i=1}^{u}(-1)^iF_i^u[p(1_{\mathcal{O}},\\
                                &\cdots,1_{\mathcal{O}},z_1,1_{\mathcal{O}},\cdots, 1_{\mathcal{O}}, z_s, 1_{\mathcal{O}},\cdots 1_{\mathcal{O}},z_{n},1_{\mathcal{O}},\cdots,1_{\mathcal{O}})]\\
                                 &= \sum\limits_{1\leq j_1 \leq\cdots \leq j_{s}\leq \cdots \leq\cdots\leq j_{n}}(-1)^{\epsilon}
                                 [\sum\limits_{i=1}^{r}(-1)^i(F_i^r p)(1_{\mathcal{O}},
                                \cdots,1_{\mathcal{O}},z_1,\\ & 1_{\mathcal{O}},\cdots, 1_{\mathcal{O}}, z_s, 1_{\mathcal{O}},\cdots
                                1_{\mathcal{O}},z_{n},1_{\mathcal{O}},\cdots, 1_{\mathcal{O}})\\&
                                 + \sum\limits_{q= 1}^{n}(-1)^{i_q} p(1_{\mathcal{O}},  \cdots, 1_{\mathcal{O}},z_1, 1_{\mathcal{O}},\cdots, 1_{\mathcal{O}},
                                (\sum\limits_{k=1}^{j_q}(-1)^{k}F_k^{j_q}(z_q)),\\ & 1_{\mathcal{O}}, \cdots, 1_{\mathcal{O}}, z_{n},1_{\mathcal{O}}, \cdots, 1_{\mathcal{O}})]\\
                                 &= (\partial_{r} p)\{z_1,\cdots, z_{n}\}+\sum\limits_{q=
                                 1}^{n}(-1)^{\delta} p\{z_1,\cdots,(\partial_{t_{q}}z_q),\cdots, z_{n}\}.
\end{aligned}
\end{equation*}
Here the sequence $(j_s)_{s= 1}^{n}$ denotes the respective
positions of the operations $(z_s)_{s= 1}^{n}$ in the right brace, $i_q$ the total number of inputs in front of $z_q$ in the right brace and $\epsilon$ has been given in the definition of the right brace structure.\\
But since by direct and arduous calculation one obtain:
$$(\partial_rp)\{z_1,\cdots,z_n\}= \left\{ \begin{array}{ll}
  0\quad \mbox{if r is even} \\
  -(F_r^r p)\{z_1,\cdots,z_n\}  \quad \mbox{else},
   \end{array} \right.$$
 then the result of second part of the theorem follows.
\end{proof}

\begin{Rem}
Notice that for a given operad $\mathcal{O}= \bigoplus_{k\geqslant
0}\mathcal{O}(k)$ with a simplicial differential right brace algebra
structure, neither the coboundary operator  $d$ nor the simplicial
differential operator $\partial$ are derivation for the
 the dot-product.
\end{Rem}

\section{Some Examples of Simplicial Operads}
 Hereafter we are going to illustrate our previous constructions by studying three distinct types of multiplicative connected operads, which
endowed with their respective right brace algebra structures are
simplicial right homotopy G-algebras as it has been shown in the
previous sections.

\subsection{Simplicial structure on the endomorphism operad $\mathcal{L}^{\prime}_{\ast}$}
For a given associative unitary $\mathbb{K}$-algebra $A$ with
multiplication $\mu_{A}$ and unit $\eta_{A}.$ The sub operad
$\mathcal{L}^{\prime}_{A}:= (\mathcal{L}^{\prime}_{A}(n))_{n\geq 0}$
 of the associated  endomorphism operad, $\mathcal{L}_{A}:= (\mathcal{L}_{A}(n))_{n\geq 0},$ is a connected unitary multiplicative
 $\Sigma$-operad. Therefore it is equipped with a simplicial structure defined as follows:\\
Set $m\circ_1 1_{0}= id_A= m\circ_2 1_{0}$ and
$D^0:=\eta:\mathcal{L}^{\prime}_{A}(0)= \mathbb{K}\longrightarrow
\mathcal{L}^{\prime}_A(1)$ such that
 $\eta(1_0)=id_A$ with $\eta$ the unit morphism of the endomorphism operad $\mathcal{L}_A.$ We have for $1\leq i\leq n,$
\begin{center}
 $F_i^n:Hom(A^{\otimes n},A)\longrightarrow(Hom(A^{\otimes (n-1)},A), n\geq 1$\\ $f\mapsto f\circ_i1_0$
\end{center}
\begin{center}
$D_i^n:Hom(A^{\otimes n},A)\longrightarrow(Hom(A^{\otimes (n+1)},A),
n\geq 1$\\ $f\mapsto f\circ_im$
\end{center}
such that for $(a_1,\cdots,a_{n-1})\in A^{\otimes(n-1)}$ and
$(a_1,\cdots,a_{n+1})\in A^{\otimes(n+1)},$ \begin{enumerate}
\item[(i)] $(f\circ_1 1_0)(a_1,\cdots,a_{n-1})=f(a_1,\cdots,a_{i-1}, 1_A,
a_{i+1},\cdots,a_{n-1})$
\item[(ii)]
$(f\circ_im)(a_1,\cdots,a_{n+1})=f(a_1,\cdots,a_{i-1},a_ia_{i+1},\cdots,a_{n+1}).$
\end{enumerate}
Let us next define some operators on the sub operad
$\mathcal{L}^{\prime}_A.$
\begin{enumerate}
\item[(a)] {\bf The dot product}
$$\begin{array}{ccc}\mathcal{L}^{\prime}_{A}\otimes\mathcal{L}^{\prime}_{A} & \stackrel{\bullet}\longrightarrow & \mathcal{L}^{\prime}_{A}\\
 f\otimes g & \mapsto & f\bullet g=(-1)^{deg(g)}m\{f,g\}
 \end{array}$$ such that for $f\in \mathcal{L}^{\prime}_{A}(r), g\in \mathcal{L}^{\prime}_{A}(s),$ $f\bullet g\in \mathcal{L}^{\prime}_{A}(r+s)$ with
$$\begin{array}{ccc}
 (f\bullet g)(a_1,\cdots,a_{r+s}) & = &
 (-1)^{rs}m(f(a_1,\cdots,a_r),g(a_{r+1},\cdots,a_{r+s}))\\
 & = & (-1)^{rs}f(a_1,\cdots,a_r)g(a_{r+1},\cdots,a_{r+s})
 \end{array}$$
\item[(b)] {\bf The coboundary operator $d$.}\\ $ For \quad all \quad n\in
\mathbb{N},$ $d_n$ is defined as follows:
$$\begin{array}{ccccc}\mathcal{L}^{\prime}_{A} & \stackrel{d_n}\longrightarrow &\mathcal{L}^{\prime}_{A}(n+1)& &,\\
f & \mapsto &
 d_n(f)& = & m\circ f+(-1)^{deg(f)}f\circ m\\
& & & = & m\circ f+(-1)^{n}\sum\limits_{i=1}^{n}(-1)^{(n-i)}f\circ_i m\\
& & & = & m\circ f + \sum\limits_{i=1}^{n}(-1)^{i}f\circ_i m\\
 & & & = & (-1)^{\mid f\mid}m\circ_1 f + m\circ_2 f + \sum\limits_{i=1}^{n}(-1)^{i}f\circ_i m
\end{array}.$$

\item[(c)] {\bf The boundary operator $\partial$}
$$for \quad all \quad n\in \mathbb{N}, \begin{array}{ccc}
\mathcal{L}^{\prime}_{A}(n)& \stackrel{\partial_n}\longrightarrow
&\mathcal{L}^{\prime}_{A}(n-1)\\
f &\mapsto &
\partial_n(f)=\sum\limits_{i=1}^{n}(-1)^if\circ_i1_0
\end{array}$$
\end{enumerate}

On the other hand, consider the Hochschild cochain complex,
$(C^*(A,A)=\mathcal{L}_A(*):=Hom(A^{\otimes *},A), d_{H}),$ of an
associative unitary $\mathbb{K}$-algebra, $A,$ with coefficients in
itself. We recall the following well known results(see \cite{S.
W.}):
\begin{enumerate}
\item[(a')] For all $n\in \mathbb{N}$, the Hochschild coboundary operator is given by\\ $d_H^{n}= \sum\limits_{i=1}^{n+1}(-1)^{i-1}\delta^{i}$ such that
 for all $f\in Hom(A^{\otimes
n},A)$,
  $$\delta^{i}(f)(a_1,\cdots,a_{n+1})=\left\{ \begin{array}{lll} a_{1} f(a_{2},\cdots,a_{n+1}) \quad \mbox{if}\quad i=1\\
   f(a_1,\cdots,a_{i-1}a_{i},
\cdots,a_{n+1})\quad \mbox{if}\quad 1<i\leq n+1 \\
f(a_1,\cdots,a_{n})a_{n+1}\quad \mbox{if}\quad i= n+2
\end{array}
\right. $$
\item[(b')] The usual cup-product is defined by: for $f\in Hom(A^{\otimes r},A), g\in Hom(A^{\otimes s},A)$,\\
 $(f\cup g)(a_1,\cdots,a_{r+s})= f(a_1,\cdots,a_r)g(a_{n+1},\cdots,a_{r+s})$
\end{enumerate}
It is obvious to check that the odot-product and coboundary operator
of the simplicial right brace algebra
$(\mathcal{L}^{\prime}_{A}(n))_{n\geq 0}$ given above coincide
respectively with usual cup-product and the external differential of
the Hochschid cochain complex(see \cite{S. W.}). Moreover the
simplicial boundary commutes with the external differential of the
Hochschid cochain complex.

\subsection{Simplicial associative operad $\mathcal{A}_{ss}.$}
\subsubsection{Partial composition of permutations}
 Let $\sum= (\sum_{n})_{n\geq 1}$ be a sequence of sets whose the $nth$ term
is the set of permutations of order $n.$ For every $1\leq i\leq n,
n\in \mathbb{N^{\ast}},$ we define on $\sum$ a partial composition,
$\begin{array}{ccc} \sum_{n} \otimes \sum_{l} &
\stackrel{\circ_i}\longrightarrow &
\sum_{n+l-1}\\
 \tau \otimes \sigma & \mapsto & \tau \circ_i \sigma
\end{array},$  in two different ways as follows:\\

\noindent{\bf 5.2.1-a- First method}\\
Set $[n]=\{1, 2,\cdots, n\}.$ Since $\tau$ and $\sigma$ are
permutations with respective orders
 $n\geq 1$ and $l\geq 1$ then
$\tau\circ_i\sigma\in\sum_{n+l-1}$ is obtained by the following
process:
\begin{enumerate}
\item[(i)]  We take the subdivision of $[n+l-1]=\{1,2, \cdots, n+l-1\}$ into n-blocks such that the ith block of length $l$ is $\{i, i+1,\cdots,
i+l-1\}.$ It is obvious that the remaining blocks are of length 1.
Thus if we denote by $B$ the above mentioned
 subdivision of $[n+l-1]$ into n-blocks, then $$B=\underbrace{(1)}_{\mbox{first block}}(2) \cdots (i-1)\underbrace{(i\quad i+1\quad \cdots \quad
i+l-1)}_{\mbox{ith block}}\underbrace{(i+1)}_{\mbox{(i + 1)th
block}}\cdots \underbrace{(n+l-1)}_{\mbox{(n-i)th block}}$$

\item[(ii)]Hereafter we permute the n blocks of $B$ by the inverse $\tau^{-1}$ of the permutation $\tau$ and next we permute
 the $ith$ block $\{i, i+1, \cdots, i+l-1\} $ by the inverse $\sigma^{-1}$ of the permutation $\sigma.$ The obtained result is denoted $\tau_{\sigma^{-1}}^{-1}(B)$.
\item[(iii)] We finally define $\tau\circ_i\sigma$ by the relation $(\tau\circ_i\sigma)^{-1}= \tau_{\sigma^{-1}}^{-1}(B).$  In other
words, for $1\leq j\leq n+l-1$, $(\tau\circ_i\sigma)(j)$ is at the
position j in $\tau_{\sigma^{-1}}^{-1}(B).$
\end{enumerate}
\textbf{Computations illustration}\\ $\tau=(4312)$, $\sigma=(231)$, $1\leq i\leq 4$, $\tau\circ_i\sigma\in\sum_6$\\
$\bullet$ Take $i=1$ and let us compute $\tau\circ_1\sigma.$\\
Since the sequence of blocks is $B=(123)(4)(5)(6),$ then
$\tau_{\sigma^{-1}}^{-1}(B)=(564\overbrace{312}).$ Thus
$\tau\circ_1\sigma=(564312)$
 where $\overbrace{312}$ is the first block of $B$ which have been permuted by $\sigma^{-1}$.\\
 If we denote graphically a permutation, $\sigma\in \sum_{n},$ by:
\begin{tikzpicture}[line cap=round,line join=round,>=triangle 45,x=1.0cm,y=1.0cm]

\clip(1.7,0.95) rectangle (6.4,2.6);

\draw (2.,1.76)-- (2.,1.);

\draw (2.,1.)-- (6.,1.);

\draw (6.,1.)-- (6.,1.82);

\draw (2.6,1.72)-- (2.58,1.);

\draw (3.1,1.44)-- (3.1,1.48);

\draw (3.86,1.46)-- (3.86,1.44);

\draw (4.78,1.48)-- (4.76,1.48);

\begin{scriptsize}

\draw[color=black] (2.04,2.25) node {$\sigma(1)$};

\draw[color=black] (4.16,-0.75) node {$g$};

\draw[color=black] (6.02,2.23) node {$\sigma(n)$};

\draw[color=black] (2.6,2.29) node {$\sigma(2)$};

\draw[color=black] (3.84,3.19) node {$j$};

\draw[color=black] (3.42,3.21) node {$k$};

\draw[color=black] (4.76,3.57) node {$l$};

\end{scriptsize}

\end{tikzpicture},
then
\begin{tikzpicture}[line cap=round,line join=round,>=triangle 45,x=1.0cm,y=1.0cm]

\clip(0.4,1.9) rectangle (9.81,3.18);

\draw (2.02,2.72)-- (2.,2.);

\draw (2.,2.)-- (4.,2.);

\draw (4.,2.)-- (4.,2.74);

\draw (3.28,2.78)-- (3.26,2.);

\draw (2.74,2.76)-- (2.74,2.);

\draw (5.02,2.76)-- (5.,2.);

\draw (5.,2.)-- (6.,2.);

\draw (6.,2.)-- (6.02,2.78);

\draw (5.54,2.76)-- (5.52,2.);

\draw (7.4,2.76)-- (7.4,2.32);

\draw (7.4,2.32)-- (7.02,2.32);

\draw (7.02,2.32)-- (7.02,2.78);

\draw (7.4,2.32)-- (7.92,2.32);

\draw (7.92,2.32)-- (7.9,2.74);

\draw (7.4,2.32)-- (7.42,1.98);

\draw (7.42,1.98)-- (9.6,2.);

\draw (9.6,2.)-- (9.58,2.76);

\draw (9.02,2.78)-- (9.000050496549402,1.9944958761151321);

\draw (8.46,2.76)-- (8.44028109745834,1.9893603770409023);

\begin{scriptsize}

\draw[color=black] (0.78,2.39) node {$\tau\circ_1\sigma$};

\draw[color=black] (2.74,3.01) node {$3$};

\draw[color=black] (3.3,2.97) node {$1$};

\draw[color=black] (2.04,2.95) node {$4$};

\draw[color=black] (1.52,2.39) node {$=$};

\draw[color=black] (4.02,2.95) node {$2$};

\draw[color=black] (5.04,2.97) node {$2$};

\draw[color=black] (5.5,2.95) node {$3$};

\draw[color=black] (4.48,2.37) node {$\circ_1$};

\draw[color=black] (6.02,3.01) node {$1$};

\draw[color=black] (7.88,4.29) node {$$};

\draw[color=black] (7.92,3.01) node {$4$};

\draw[color=black] (6.52,2.35) node {$=$};

\draw[color=black] (8.44,4.33) node {$$};

\draw[color=black] (6.98,3.01) node {$5$};

\draw[color=black] (7.42,2.99) node {$6$};

\draw[color=black] (8.48,2.97) node {$3$};

\draw[color=black] (9.58,3.01) node {$2$};

\draw[color=black] (9.02,3.01) node {$1$};

\end{scriptsize}

\end{tikzpicture}\\
$\bullet$ Take now $i=2$ and let us compute $\tau\circ_2\sigma.$\\
 The sequence of blocks is $B=(1)(234)(5)(6)$ then $\tau_{\sigma^{-1}}^{-1}(B)=(56\overbrace{423}1).$
Thus $\tau\circ_2\sigma=(645312)$
with $\overbrace{423}$ the second block of $B$ which have been permuted by $\sigma^{-1}$. Hence we obtain graphically \\
\begin{tikzpicture}[line cap=round,line join=round,>=triangle 45,x=1.0cm,y=1.0cm]

\clip(0.4,1.95) rectangle (9.8,3.17);

\draw (2.02,2.72)-- (2.,2.);

\draw (2.,2.)-- (4.,2.);

\draw (4.,2.)-- (4.,2.74);

\draw (3.28,2.78)-- (3.26,2.);

\draw (2.74,2.76)-- (2.74,2.);

\draw (5.02,2.76)-- (5.,2.);

\draw (5.,2.)-- (6.,2.);

\draw (6.,2.)-- (6.02,2.78);

\draw (5.54,2.76)-- (5.52,2.);

\draw (7.02,2.76)-- (7.,2.);

\draw (7.,2.)-- (8.,2.);

\draw (8.,2.)-- (8.04,2.32);

\draw (7.48,2.74)-- (7.48,2.32);

\draw (7.48,2.32)-- (8.04,2.32);

\draw (8.02,2.74)-- (8.04,2.32);

\draw (8.04,2.32)-- (8.5,2.3);

\draw (8.48,2.72)-- (8.5,2.3);

\draw (8.,2.)-- (9.56,2.);

\draw (9.56,2.)-- (9.58,2.74);

\draw (9.02,2.72)-- (9.,2.);

\begin{scriptsize}

\draw[color=black] (0.78,2.39) node {$\tau\circ_2\sigma$};

\draw[color=black] (2.72,3.01) node {$3$};

\draw[color=black] (3.26,2.97) node {$1$};

\draw[color=black] (2.04,2.97) node {$4$};

\draw[color=black] (1.52,2.39) node {$=$};

\draw[color=black] (4.04,2.99) node {$2$};

\draw[color=black] (4.56,2.35) node {$\circ_2$};

\draw[color=black] (5.54,2.99) node {$3$};

\draw[color=black] (5.02,2.99) node {$2$};

\draw[color=black] (6.02,3.01) node {$1$};

\draw[color=black] (6.56,2.35) node {$=$};

\draw[color=black] (8.5,0.97) node {$$};

\draw[color=black] (7.,3.01) node {$6$};

\draw[color=black] (7.9,0.69) node {$$};

\draw[color=black] (7.46,2.99) node {$4$};

\draw[color=black] (8.02,2.97) node {$5$};

\draw[color=black] (8.44,2.99) node {$3$};

\draw[color=black] (8.98,2.97) node {$1$};

\draw[color=black] (9.58,3.03) node {$2$};

\draw[color=black] (9.11,1.04) node {$$};

\end{scriptsize}

\end{tikzpicture}\\ This means that $\sigma$ has been  inserted at the position 2 of the representation of $\tau$.\\

\noindent{\bf 5.2.1-b-Second method}\\ Once more let $\tau$ and
$\sigma$ be two permutations with respective orders $n$ and $l.$ For
all $1\leq i\leq n,$  $\tau\circ_i\sigma\in\sum_{n+l-1}$ can be also
computed as follows:
\begin{enumerate}
\item[(i)] We start by computing $\tau(i)$ then  consider the block $A=(\tau(i),\quad \tau(i)+1, \quad ...\quad, \tau(i)+l-1)$ with length $l$
that is permuted by $\sigma$ and the result is denoted $\sigma(A)$.
\item[(ii)] Hereafter we consider $B,$ a subdivision of  $[n+l-1]$ into n-blocks where the block $\sigma(A)$ is located at the position $\tau(i)$.
\item[(iii)] We end by permuting by $\tau$ denoted here by $\tau^{B},$ the previous result (ii), that is we permute the n-block of $B$ as components
 of the permutation by $\tau;$ and therefore the obtained result is $\tau\circ_i\sigma$.
\end{enumerate}
This second method leads us to the following formal definition of
the partial composition of permutations:

$$ \begin{array}{c}
\circ_i:\sum_n\times\sum_l  \longrightarrow  \sum_{n+l-1}, 1\leq i\leq n\\
 (\tau,\sigma) \mapsto  \tau\circ_i\sigma,\\
  \quad \mbox{such that for all} \quad  1\leq j\leq n+l-1,
 \end{array}$$

$$\tau\circ_i\sigma(j)=\left\{\begin{array}{lll}
\tau^B(j)\quad if \quad 1\leq j\leq i-1\\
\tau(i)-1 + \sigma(j-i+1)\quad if \quad i\leq j\leq i+l-1\\
\tau^B(j-l+1)\quad if\quad i+l\leq j\leq n+l-1
\end{array}
\right.$$

 \textbf{Computations illustration.}\\ Let $\tau=(4312)\in \sum_{4}$ and
$\sigma=(231)\in \sum_{3}$, then we are going to compute
consecutively $\tau\circ_1\sigma\in \sum_{6}$ and $\tau\circ_2\sigma\in\sum_6$\\
 $\bullet$ Since $\tau(1)=4$ and $A=(456)$ then $\sigma(A)=(564)$. Thus the subdivision of $[6]$ into 4 blocks is $B=(1)(2)(3)\sigma(A)=(1)(2)(3)(564)$
  and therefore $\tau(B)=(564312)=\tau\circ_1\sigma$.\\
$\bullet$ On other hand $\tau(2)=3$  and $A=(345)$ implies that
$\sigma(A)=(453)$. Now considering the subdivision of $[6]$ into 4
blocks $B=(1)(2)\sigma(A)(6)=(1)(2)(453)(6),$ we obtain
$\tau(B)=(645312)=\tau\circ_2\sigma$.

 \subsubsection{ Construction of the Simplicial Structure on the unital associative Operad}
 For a permutation $\tau\in\sum_n$, we denote $\tau$ by $(\tau(1),\cdots,\tau(n))$.
 \begin{De}
 If u is a sequence of integers without multiple occurrence of the same integer,
 the standardization std(u) of u is the unique permutation $\sigma\in\sum_l$ order-isomorphic to u.\\ For instance $std(291847)=(261534)\in\sum_6$
 \end{De}
  Let $\mathcal{A}_{ss}= \bigoplus_{n\geq 0}\mathbb{K}[\sum_n]$ denotes the unitary multiplicative and
connected associative operad with
$1_{\mathcal{A}_{ss}}\in\mathcal{A}_{ss}(1)\cong \mathbb{K},$ its
unit, $1_{\mathbb{K}}\in \mathbb{K}\cong\mathcal{A}_{ss}(0)$ and
$m=(12)\in\mathbb{K}[\sum_2]=\mathcal{A}_{ss}(2),$ its
multiplication. It is obvious that $m$ satisfies
the relation: $m\circ_1 1_{\mathbb{K}}= 1_{\mathcal{A}_{ss}}= m\circ_2 1_{\mathbb{K}}$.\\
Hence for $n\in \mathbb{N}^{\ast}$ and $1\leq i\leq n,$ we define
respectively the face morphism, $F_{i}^{n},$ and the degeneracy
morphism, $D_{i}^{n},$ as
 follows:\\
$\bullet$
 $\begin{array}{ccc}
 \mathbb{K}[\sum_n] & \stackrel{F_{i}^{n}}\longrightarrow & \mathbb{K}[\sum_{n-1}],\\
\tau & \longmapsto & \tau\circ_i1_{\mathbb{K}}
\end{array}$
where for $1\leq i\leq n$, $\tau\circ_i1_{\mathbb{K}}$ is the
standardized of
$(\tau(1),\cdots,\tau(i-1),\tau(i+1),\cdots,\tau(n))$ denoted by $std(\tau(1),\cdots,\tau(i-1),\tau(i+1),\cdots,\tau(n))$.\\

$\bullet$ $\begin{array}{ccc} \mathbb{K}[\sum_n]& \stackrel{D_i^n}\longrightarrow & \mathbb{K}[\sum_{n+1}]\\
\tau & \longmapsto & \tau\circ_i m
\end{array}$\\
Observe moreover that if one denotes by $ D^0:=\eta,$ the unit
homomorphism of the operad $\mathcal{A}_{ss},$ then
$\eta(1_{\mathbb{K}})= D^0(1_{\mathbb{K}})= 1_{\mathcal{A}_{ss}}$
and $D^1(1_{\mathcal{A}_{ss}})= 1_{\mathcal{A}_{ss}}\circ_1 m = m.$

\textbf{Computations illustration.}\\
 Let $\tau =(4312)\in  \mathbb{K}[\sum_4],$ we have:\\

\begin{tikzpicture}[line cap=round,line join=round,>=triangle 45,x=1.0cm,y=1.0cm]

\clip(1.7,0.9) rectangle (11.9,8.2);

\draw (4.02,7.76)-- (4.,7.);

\draw (4.,7.)-- (5.58,7.02);

\draw (5.58,7.02)-- (5.6,7.74);

\draw (4.99983979493752,7.012656199935917)-- (5.02,7.74);

\draw (4.51991669336751,7.006581223966677)-- (4.52,7.74);

\draw (4.,5.76)-- (4.,5.);

\draw (4.,5.)-- (5.58,5.);

\draw (5.58,5.)-- (5.58,5.76);

\draw (5.02,5.76)-- (5.,5.);

\draw (4.5,5.74)-- (4.52,5.);

\draw (4.02,3.78)-- (4.,3.);

\draw (4.,3.)-- (5.58,3.);

\draw (5.58,3.)-- (5.6,3.76);

\draw (5.02,3.74)-- (5.,3.);

\draw (4.5,3.)-- (4.52,3.78);

\draw (4.02,1.76)-- (4.,1.);

\draw (4.,1.)-- (5.56,1.);

\draw (5.56,1.)-- (5.56,1.74);

\draw (5.02,1.76)-- (5.,1.);

\draw (4.54,1.76)-- (4.52,1.);

\draw (4.51991669336751,7.006581223966677)--
(4.519916693367512,7.006581223966679);

\draw (4.8180967638577386,7.010355655238706)-- (4.,7.);

\draw (5.6,7.74)-- (5.591642251349267,7.439121048573632);

\draw (5.58,5.76)-- (5.58,5.16);

\draw (8.58,7.74)-- (8.58,7.);

\draw (8.58,7.)-- (9.5,7.);

\draw (9.5,7.)-- (9.5,7.76);

\draw (8.58,7.74)-- (8.58,7.);

\draw (9.,7.)-- (9.02,7.74);

\draw (9.02,7.74)-- (9.,7.);

\draw (5.58,7.02)-- (5.593862760215883,7.5190593677717805);

\draw (5.58,5.)-- (5.58,5.48);

\draw (5.58,3.)-- (5.6,3.76);

\draw (5.56,1.)-- (5.56,1.74);

\draw (8.54,5.74)-- (8.54,5.);

\draw (8.54,5.)-- (9.36,5.);

\draw (9.36,5.)-- (9.36,5.74);

\draw (9.02,5.76)-- (9.,5.);

\draw (9.,5.)-- (9.36,5.);

\draw (9.,5.)-- (8.54,5.);

\draw (8.54,5.)-- (8.54,5.74);

\draw (9.36,5.)-- (9.36,5.74);

\draw (10.78,5.74)-- (10.78,4.98);

\draw (10.78,4.98)-- (11.66,5.);

\draw (11.66,5.)-- (11.68,5.76);

\draw (11.221135776974704,4.9900258131130615)-- (11.24,5.76);

\draw (8.54,3.76)-- (8.54,3.);

\draw (8.54,3.)-- (9.32,3.);

\draw (9.32,3.)-- (9.32,3.74);

\draw (9.,3.72)-- (9.,3.);

\draw (10.8,3.78)-- (10.8,2.98);

\draw (10.8,2.98)-- (11.64,2.98);

\draw (11.3,3.72)-- (11.3,2.98);

\draw (11.64,2.98)-- (11.64,3.74);

\draw (8.44,1.78)-- (8.44,1.);

\draw (8.44,1.)-- (9.32,1.);

\draw (9.32,1.)-- (9.32,1.76);

\draw (8.92,1.76)-- (9.,1.);

\draw (9.,1.)-- (8.44,1.);

\draw (10.78,1.76)-- (10.8,1.);

\draw (10.8,1.)-- (11.66,1.);

\draw (11.66,1.)-- (11.66,1.8);

\draw (11.2,1.78)-- (11.22,1.);

\begin{scriptsize}

\draw[color=black] (2.18,7.61) node {$\tau\circ_11_{\mathbb{K}}$};

\draw[color=black] (5.,6.11) node {$1$};

\draw[color=black] (5.,8.09) node {$1$};

\draw[color=black] (4.52,8.11) node {$3$};

\draw[color=black] (4.06,8.11) node {$4$};

\draw[color=black] (2.14,5.53) node {$\tau\circ_21_{\mathbb{K}}$};

\draw[color=black] (5.04,4.11) node {$1$};

\draw[color=black] (5.54,4.03) node {$2$};

\draw[color=black] (3.98,6.09) node {$4$};

\draw[color=black] (3.54,5.53) node {$=$};

\draw[color=black] (2.16,3.59) node {$\tau\circ_31_{\mathbb{K}}$};

\draw[color=black] (5.6,2.11) node {$2$};

\draw[color=black] (4.48,4.11) node {$3$};

\draw[color=black] (3.48,3.53) node {$=$};

\draw[color=black] (3.98,4.13) node {$4$};

\draw[color=black] (2.14,1.51) node {$\tau\circ_41_{\mathbb{K}}$};

\draw[color=black] (4.02,2.15) node {$4$};

\draw[color=black] (5.,2.15) node {$1$};

\draw[color=black] (4.54,2.13) node {$3$};

\draw[color=black] (3.51,1.52) node {$=$};

\draw[color=black] (4.53,6.2) node {$3$};

\draw[color=black] (3.47,7.64) node {$=$};

\draw[color=black] (5.69,8.14) node {$2$};

\draw[color=black] (5.65,6.16) node {$2$};

\draw[color=black] (7.01,7.52) node {$1_{\mathbb{K}}$};

\draw[color=black] (9.09,8.16) node {$1$};

\draw[color=black] (9.51,8.16) node {$2$};

\draw[color=black] (8.11,7.54) node {$=$};

\draw[color=black] (8.57,8.16) node {$3$};

\draw[color=black] (6.01,7.5) node {$\circ_1$};

\draw[color=black] (6.01,5.56) node {$\circ_2$};

\draw[color=black] (6.05,3.54) node {$\circ_3$};

\draw[color=black] (6.03,1.5) node {$\circ_4$};

\draw[color=black] (7.05,5.54) node {$1_{\mathbb{K}}$};

\draw[color=black] (9.41,6.18) node {$2$};

\draw[color=black] (10.81,6.14) node {$3$};

\draw[color=black] (8.11,5.54) node {$=$};

\draw[color=black] (10.07,5.54) node {$=$};

\draw[color=black] (8.55,6.14) node {$4$};

\draw[color=black] (9.03,6.16) node {$1$};

\draw[color=black] (11.27,6.14) node {$1$};

\draw[color=black] (11.77,6.14) node {$2$};

\draw[color=black] (10.81,4.18) node {$3$};

\draw[color=black] (11.71,4.12) node {$1$};

\draw[color=black] (11.27,4.14) node {$2$};

\draw[color=black] (7.05,3.56) node {$1_{\mathbb{K}}$};

\draw[color=black] (8.57,4.18) node {$4$};

\draw[color=black] (9.01,4.14) node {$3$};

\draw[color=black] (8.07,3.56) node {$=$};

\draw[color=black] (9.39,4.12) node {$2$};

\draw[color=black] (11.31,2.88) node {$$};

\draw[color=black] (10.05,3.6) node {$=$};

\draw[color=black] (12.05,3.58) node {$$};

\draw[color=black] (7.01,1.58) node {$1_{\mathbb{K}}$};

\draw[color=black] (8.95,2.12) node {$3$};

\draw[color=black] (9.35,2.16) node {$1$};

\draw[color=black] (8.07,1.5) node {$=$};

\draw[color=black] (8.47,2.2) node {$4$};

\draw[color=black] (10.03,1.58) node {$=$};

\draw[color=black] (11.25,2.16) node {$2$};

\draw[color=black] (11.71,2.16) node {$1$};

\draw[color=black] (10.83,2.12) node {$3$};

\end{scriptsize}
\end{tikzpicture}\\
\textbf{Computations without using graphic:} let $\tau=(4312)\in S_4$, for $1\leq i\leq 4$ we have\\
 $\tau\circ_11_{\mathbb{K}}=\tau\mid_{\{2<3<4\}}=std(312)=(312)$\\
 $\tau\circ_21_{\mathbb{K}}=\tau\mid_{\{1<3<4\}}=std(412)=(312)$\\
 $\tau\circ_31_{\mathbb{K}}=\tau\mid_{\{1<2<4\}}=std(432)=(321)$\\
 $\tau\circ_41_{\mathbb{K}}=\tau\mid_{\{1<2<3\}}=std(431)=(321)$\\

The above constructed face and degeneracy maps $F_i^n$ and $D_i^n$
satisfy the following properties:
\begin{center}
$F_iF_j=F_{j-1}F_i$, $i<j$ \\
$D_iD_j=D_{j+1}D_i$, $i\leq j$\\
$F_iD_j=\left\lbrace
\begin{array}{ll}
D_{j-1}F_i, i<j\\ I, i\in\{j,j+1\}\\ D_jF_{i-1},i>j+1
\end{array}
\right.$
\end{center}
Indeed these properties follow from associativity on the operad
$\mathcal{A}_{ss}$ hence the associative operad $\mathcal{A}_{ss}$
is a simplicial operad.\\ Moreover for two generators $\tau$ and
$\sigma$ of $\mathcal{A}_{ss}$, we have
\begin{center}
$\tau\bullet\sigma=(-1)^{(\mid\tau\mid+1)(\mid\sigma\mid+1)}\tau\odot\sigma=(-1)^{(\mid\tau\mid+1)(\mid\sigma\mid+1)}\tau\star\sigma$
\end{center}
where $\star$ represent the concatenation of $\tau$ and $\sigma.$ So
the odot product $\odot$ coincides with $\star$.
\subsubsection{Analogy between the Malvenuto-Reutenauer coproduct, $\Delta_{MR}$, and $\Delta_{\mathcal{A}_{ss}}$.}
In 1995, Malvenuto and Reutenauer defined on $\mathcal{A}_{ss}=
\bigoplus\limits_{n\in\mathbb{N}}\mathbb{K}[\Sigma_n]$ a Hopf
algebra structure called shuffle Hopf algebra of permutations or
malvenuto-Reutenauer Hopf algebra(see \cite{M.R}). Having shown that
an arbitrary connected and multiplicative operad has a graded
differential coalgebra structure, we show in particular that the
coalgebra structure on the unit associative operad
$\mathcal{A}_{ss}=\bigoplus\limits_{n\in\mathbb{N}}\mathbb{K}[\Sigma_n]$
coincides with the existing coalgebra structure on the
Malvenuto-Reutenauer Hopf algebra.\\
Let us first recall the coalgebra structure defined by  Malvenuto
and Reutenauer.

\begin{De}
    The standardization map is the map associating to a word of length k, $X=x_1\cdots x_k$ in the letters $1,\cdots,k$,
    the unique word $st(X)=y_1\cdots y_k$ in the letters $1,\cdots,k$, without repetition of the letters, such that the
    relative order of the letters is preserved:$y_i<y_j$ if $x_i<x_j$ and such that furthermore, if $x_i=x_j$ and
    $i<j,$
    then $y_i<y_j$.
\end{De}
\begin{Ex} st(3745)=(1423), st(2221)=(2341).
 \end{Ex}
\indent Heather we are going to recall properly the
Malvenuto-Reutenauer coproduct. For this, consider
$\sigma=(\sigma(1),\cdots,\sigma(n))\in \Sigma_n$, we write
$\sigma_i^-$ for the word $\sigma(1)\cdots\sigma(i)$ and
$\sigma_{i+1}^+$ for the word $\sigma(i+1)\cdots\sigma(n)$.\\ The
coassociative coproduct of Malvenuto-Reutenauer Hopf algebra is
defined as follows: for all $\sigma\in \Sigma_n$
\begin{center}
    $\Delta_{MR}=\sum\limits_{i=0}^nst(\sigma_i^-)\otimes st(\sigma_{i+1}^+)$.
\end{center}
For example let $\sigma=(3124)\in \Sigma_4$
\begin{align*}
    \Delta_{MR}(\sigma)&=1_{\mathbb{K}}\otimes st(3124)+st(3)\otimes st(124)+st(31)\otimes st(24)\\
    &+st(312)\otimes st(4)+st(3124)\otimes 1_{\mathbb{K}};\\
    &=1_{\mathbb{K}}\otimes(3124)+(1)\otimes(123)+(21)\otimes(12)+(312)\otimes(1)+(3124)\otimes 1_{\mathbb{K}}
\end{align*}
And the counit map defined as
\begin{center}
    $\varepsilon:\mathcal{A}_{ss} \longrightarrow\mathbb{K}$
\end{center}
such that $\varepsilon_0=Id$ and $\varepsilon_n=0, \forall n\neq
0$.\\ Then $(\mathcal{A}_{ss},\Delta_{MR},\varepsilon)$ is the
coalgebra in the shuffle Hopf algebra of permutations.\\

Now coming back to the differential graded coalgebra structure on
the connected multiplicative operad $\mathcal{A}_{ss}$,
$(\mathcal{A}_{ss},\Delta_{\mathcal{A}_{ss}},\varepsilon,\partial),$
 where $\forall \sigma\in \Sigma_n$,
$\Delta_{\mathcal{A}_{ss}}(\sigma)\sum\limits_{j=0}^{n}\widetilde{F}^{n-j}_{n}\sigma\otimes\widetilde{F}_1^j\sigma.$
Here $\widetilde{F}^{n-j}_{n}=F_{j+1}\cdots F_n$,
$\widetilde{F}_1^j=\underbrace{F_1\cdots F_1}_{j-times}$,
$\widetilde{F}^{0}_{n}=0=\widetilde{F}_1^0$
    with
$$\begin{array}{ccc}
    \mathbb{K}[\Sigma_n]&\stackrel{F_{i}}\longrightarrow &\mathbb{K}[\Sigma_{n-1}], 1\le i\le n\\
    \sigma& \mapsto &
    F_i\sigma=\sigma\mid_{[n]-\{i\}}=st(\sigma(1),\cdots,\sigma(i-1),\sigma(i+1),\cdots,\sigma(n)).
\end{array}$$
We have that $\Delta_{MR}$ and $\Delta_{\mathcal{A}_{ss}}$ coincide.
 Indeed for all $\sigma\in \Sigma_n$, $1\le j\le n$,
\begin{align*}
    \widetilde{F}^{n-j}_{n}\sigma &=F_{j+1}\cdots F_n\sigma;\\
    &=F_{j+1}\cdots F_{n-1}(st(\sigma(1),\cdots,\sigma(n-1)));\\
    &=\underbrace{st(\cdots(st(\sigma(1),\cdots,\sigma(n-1)))\cdots)}_{(n-j-1)-times};\\
    &=st(\sigma(1),\cdots,\sigma(j));\\
    &=st(\sigma_j^-)
\end{align*}
and
\begin{align*}
    \widetilde{F}_1^j\sigma &=\underbrace{F_1\cdots F_1}_{j-times}\sigma;\\
    &=\underbrace{F_1\cdots F_1}_{(j-1)-times}(st(\sigma(2),\cdots,\sigma(n)));\\
    &=\underbrace{st(\cdots(st(\sigma(2),\cdots,\sigma(n)))\cdots)}_{(n-j-1)-times};\\
    &=st(\sigma(j+1),\cdots,\sigma(n));\\
    &=st(\sigma_{j+1}^+)
\end{align*}

\subsection{Simplicial Shift operads}
We show through the following statements that the connectedness
condition used in our construction of simplicial structure on
operads is not always necessary.
\subsubsection{Construction of a shift operation on $\mathbb{N}^{\ast}.$}
Denote by $\mathcal{E}=\{E_n\}_{n\geq 0}$ the family of sets such
that for all $n\geq 1$, $E_n$ is the set of subsets of
$\mathbb{N}^{\ast}$ with n distinct  ordered elements with respect
to the relation $\le$, that is $E_n=\{X_{n} = \{a_1<\cdots\le a_n\}
\mid a_i\in \mathbb{N}^{\ast}, a_i\neq a_{i+1}, 1\leq i\leq n\}$ and
since the empty set is a subset of $\mathbb{N}^{\ast}$ with zero
element, we set by convention: $E_0=\{X_0= \emptyset= 0\}$.
\begin{De}
 Let $n\geq 1$ and  $E_n,$  the set of subsets of $\mathbb{N}^{\ast}$ with n
distinct  ordered elements with respect to the relation $\le.$ The
shift-operation on $E_n$ denoted $\oplus$ is the operation:
$$\begin{array}{ccccc}
\oplus :& E_n\times\mathbb{Z} & \longrightarrow & E_n &\\
&(X_n,p)& \mapsto & X_n\oplus p= & a_1+p\le\cdots\le a_n+p
\end{array}$$
with the condition: $X_n\oplus p=X_n$ if there exists $1\leq i\leq
n$ such that $a_i+p\leq 0.$
\end{De}
\begin{Rem}
\begin{enumerate}
 \item With the convention $X_0\in E_0, X_0+p=X_0,$ $\mathcal{E}=\{E_n\}_{n\geq 0}$ is invariant under the shift-operation.
 \item For all $n\in\mathbb{N}^{\ast},$  there exists on $E_n$ the structure of poset given the relation $\mathcal{R}$ defined as follows: for all
 $X_n=\{a_1\le\cdots\le a_n\}$ and $Y_n=\{b_1\le\cdots\le b_n\}$,
$X_n\mathcal{R} Y_n$ if and only if  $a_i\le b_i,$ for all $1\le
i\le n.$
\end{enumerate}
\end{Rem}
\noindent{\bf Notations} - Let $X_n=\{a_1\le\cdots\le a_n\}\in E_n$
be a subset of  $\mathbb{N}^{\ast}$ with n distinct ordered
elements. $X_n$ is simply denoted
 $X_n=\{a_1\le\cdots\le a_n\}$.\\
- For all $X_n=(a_1\le\cdots\le a_n)\in E_n$, $(X_n^{\hat{a}_i}\in
E_{n-1})_{ 1\leq i\leq n}$ denotes the family  of elements obtained
from $X_n$ by deleting $a_i$, that is
$X_n^{\hat{a}_i}=(a_1\le\cdots\le
a_{i-1}\le a_{i+1}\le\cdots\le a_n)\in E_{n-1}$\\
 More generally for all $2\le k\le n,$
$X_n^{\hat{a}_{i_1},\hat{a}_{i_2},\cdots \hat{a}_{i_k}}=(a_1\le
..\le a_{i_{1}-1}\le a_{i_{1}+1}\le\cdots\le a_{i_{2}-1}\le
a_{i_{2}+1}\le ..\le a_{i_{k}-1}\le a_{i_{k}+1}\le ..\le a_n)\in
E_{n-k}$.
\subsubsection{Operadic structure on $\mathcal{E}=\{E_n\}_{n\geq 0},$} In the sequel we are
going to define on $\mathcal{E}=\{E_n\}_{n\geq 0},$ an operadic
structure called shift-operad of finite ordered subset of
$\mathbb{N}^{\ast}$ as follows:\\ \indent{\bf 5.3.2.a- Partial
composition $\circ_i$ on $\mathcal{E}$}\\ Let $p, q \in
\mathbb{N}^{\ast},$ the partial composition of $E_p$ and $E_q$
consists to  insert into an element $X_p=(a_1\leq ..\le a_p)\in E_p$
at  a given position $i$, an element  $Y_q=(b_1\le ..\le b_q)\in
E_q$ by applying the shift-operation and the resulted element is
denoted $X_p\circ_i Y_q \in E_{p+q-1}$. More precisely for $1\le
i\le p,$ we have:
$$\begin{array}{ccccc}
 \circ_i :& E_p\times E_q & \longrightarrow & E_{p+q-1}&  \\
& (X_p,Y_q) & \mapsto & X_p\circ_i Y_q=&
X_p^{\hat{a}_i,...,\hat{a}_p}\le Y_q\oplus(a_i-1)\le
X_p^{\hat{a}_1,...,\hat{a}_i}\oplus (b_q-1).
\end{array}$$
That is\\
$X_p\circ_i Y_q=(\underbrace{a_1\le ..\le
a_{i-1}}_{X_p^{\hat{a}_i,...,\hat{a}_p}}\le\underbrace{b_1+a_i-1\le
..\le b_q+a_i-1}_{Y_q\oplus (a_i-1)}
\le\underbrace{a_{i+1}+b_q-1\le ..\le a_p+b_q-1}_{X_p^{\hat{a}_1,...,\hat{a}_i}\oplus(b_q-1)})$.\\
\indent{\bf 5.3.2.b-composition product on $\mathcal{E}$.}\\ The
composition product on $\mathcal{E}=\{E_n\}_{n\geq 0}$ is defined by
$$\begin{array}{cccc}
\gamma:& E_n\times E_{t_1}\times\cdots\times E_{t_n}& \longrightarrow & E_{t_1+t_2+\cdots +t_n}\\
&(X_n, X_{t_1},\cdots,X_{t_n})&\mapsto & \gamma(X_n,
X_{t_1},\cdots,X_{t_n}).
\end{array}$$
such that if $X_n=(a_1\le\cdots\le a_n),$ for all $1\le i\le n,$
$X_{t_i}=(a_1^i\le\cdots\le a_{t_i}^i)$ we have
$$\begin{array}{ccc}
\gamma(X_n, X_{t_1},\cdots,X_{t_n})& = & X_{t_1}\oplus(a_1-1)\le
X_{t_2}\oplus(a_2+a_{t_1}^1-2)\le ..\le
X_{t_i}\oplus\\
& & (a_i+\sum\limits_{j=1}^{i-1}a_{t_j}^j-i)\le ..\le
X_{t_n}\oplus(a_n+\sum\limits_{j=1}^{n-1}a_{t_j}^j-n)
\end{array}$$
\begin{Ex}
Let $X_3=(a_1\le a_2\le a_3), X_{t_1}=(b_1\le ...\le b_5), X_{t_2}=
(c_1\le ...\le c_4), X_{t_3}=(d_1\le d_2).$ One has
$$\begin{array}{ccc}
\gamma(X;X_{t_1},X_{t_2},X_{t_3})& = & X_{t_1}\oplus(a_1-1)\le X_{t_2}\oplus(a_2+b_5-2)\le X_{t_3}\oplus(a_3+b_5+c_4-3)\\
& = & (\underbrace{b_1+a_1-1\le ..\le
b_5+a_1-1}_{X_{t_1}\oplus(a_1-1)}\le\underbrace{c_1+a_2+b_5-2\le
..\le c_4+a_2+b_5-2}_{X_{t_2}\oplus(a_2+b_5-2)}\le\\ & &
\underbrace{d_1+a_3+b_5+c_4-3\le ..\le
d_2+a_3+b_5+c_4-3}_{X_{t_3}\oplus(a_3+b_5+c_4-3)})
\end{array}$$
\end{Ex}

\begin{Rem}
\begin{enumerate}
\item Set $E_1=\{\{i\} \mid i\in\mathbb{N}^*\}.$ $E_1$ is equipotent to  $\mathbb{N}^{\ast}$ and the unit element of $\mathcal{E}$
denoted $\mathbf{1}_E\in E_1$ is the subset $\mathbf{1}_E= \{1\}.$
Indeed for all $n\in \mathbb{N}^{\ast},$ consider an arbitrary
$X_n\in E_n$ and $1\le i\le n.$ We have $X_n \circ_i \mathbf{1}_E =
X_n= \mathbf{1}_E \circ_1 X_n$.
\item The composition product and the partial composition are linked
as follows: let $X\in E_n$ and $Y\in E_m.$ One has for $1\le i\le
n$,
 $X\circ_i Y=
 \gamma(X;\underbrace{\mathbf{1}_E,\cdots,\mathbf{1}_E}_{(i-1)-times},Y,\underbrace{\mathbf{1}_E,\cdots,\mathbf{1}_E}_{(n-i)-times}).$\\
More precisely if $X=(a_1\le\cdots\le a_n), \quad Y=(b_1\le\cdots\le
b_m),$ $$\begin{array}{ccc}
\gamma(X;\underbrace{\mathbf{1}_E,\cdots,\mathbf{1}_E}_{(i-1)-times},Y,\underbrace{\mathbf{1}_E,\cdots,\mathbf{1}_E}_{(n-i)-times})&
 = & (a_1\le\cdots\le a_{i-1}\le Y\oplus(a_i+\sum\limits_{j=1}^{i-1}1-i)\le
  \mathbf{1}_E\oplus(a_{i+1}+\\ & &
  \sum\limits_{j=1}^{i-1}1-i-1+b_m)\le\cdots\le\mathbf{1}_E\oplus(a_n+b_m-2))\\&&\\
& = &(a_1\le\cdots\le a_{i-1}\le Y\oplus(a_i-1)\le \\ && \\& &
a_{i+1}+b_m-1\le\cdots\le a_n+b_m-1)\\&&\\ & = &  X\circ_i Y
\end{array}$$

\item Since for all $n\in \mathbb{N}^{\ast},$ $X_n\in E_n$ is an ordered subset the symmetric group
 $\Sigma_n$ does not act on $E_n$. So the operad $\mathcal{E}=\{E_n\}_{n\geq
0}$ is a non symmetric operad.
\item Consider $M\in E_2$ such that $M=1\le 2$. $M$ satisfies the
following identities: $M\circ_1 M = M\circ_2 M$. Indeed
$M\circ_1M=1\le 2\le 3=M\circ_2M$. Thus $\mathcal{E}=\{E_n\}_{n\geq
0},$ is multiplicative operad with multiplication $M.$
\end{enumerate}
\end{Rem}

\begin{Pro}
The composition product given above defines a non symmetric
multiplicative operad called {\bf shift-operad} of finite ordered
subsets of $\mathbb{N}^{\ast}$
\end{Pro}
\begin{proof}
To show that the composition product is associative it is enough
prove that the partial composition is associative.\\
So let $X_n=(a_1\le\cdots\le a_n)\in E_n, Y_m=(b_1\le\cdots\le
b_m)\in E_m,$ and $Z_l=(c_1\le\cdots\le c_l)\in E_l.$ We have
$$\begin{array}{ccc}
 X\circ_i(Y\circ_j Z)&= &  X\circ_i(b_1 \le \cdots \le b_{j-1}
\le b_j + c_1 - 1 \le \cdots \\  & & \le b_j + c_l - 1\le
b_{j+1}+c_l-1\le \cdots \le b_m + c_l - 1) \\  &= &  a_1\le\cdots\le
a_{i-1}\le (Y\circ_jZ)\oplus(a_i-1)\le a_{i+1} \\  & + &
b_m+c_l-2\le\cdots\le a_n+b_m+c_l-2 \\ & = &
a_1\le\cdots\le a_{i-1}\le b_1+a_i-1\le\cdots\le b_{j-1}+a_i-1\le b_j\\
 & + & c_1+a_i-2\le\cdots\le b_j+c_l+a_i-2\le b_{j+1}+c_l+a_i-2\le\cdots\le b_m \\ & + & c_l+a_i-2\le a_{i+1}+b_m+c_l-2\le\cdots\le a_n+b_m+c_l-2 \\
&  = &  (X\circ_i Y)\circ_{j+i-1}Z
\end{array}$$

On the other hand if $i <j,$ then
$$\begin{array}{ccc} (X\circ_jZ)\circ_i Y &= & (a_1\le\cdots\le
a_{j-1}\le c_1+a_j-1\le\cdots  \le c_l+a_j-1\le \\ & & a_{j+1} +
c_l-1\le\cdots\le a_n + c_l-1)\circ_iY\\ && \\ & = & a_1\le\cdots\le
a_{i-1}\leq b_1+a_i-1\le\cdots \le b_m+a_i-1\le\\
&  & a_{i+1}+b_m-1\le\cdots\le a_{j-1}+b_m-1\le c_1+a_j+b_m-2\le\cdots\\
& &\le c_l+a_j + b_m-2\le a_{j+1}+c_l+b_m-2\le\cdots\le a_n +
c_l+b_m-2\\ && \\ & = & (X\circ_iY)\circ_{j+m-1}Z
\end{array}$$
Hence the partial composition is associative and therefore so is the
composition product.
\end{proof}
\subsubsection{Simplicial structure on $\mathcal{E}= \{E_n\}_{n\geq
0}.$}
 Since the operad $\mathcal{E}=\{E_n\}$ has a multiplication, we are
 going to equip it with a simplicial structure as follows
 \begin{Pro}
Let $\mathcal{E}= \{E_n\}_{n\geq 0}$ be endowed with the operad
structure given above. $\mathcal{}E= \{E_n\}_{n\geq 0}$ has a
simplicial
structure whose face and degeneracy morphisms are:\\
For $1\le i\leq n,$
$$\begin{array}{cccc}
D_i :& E_n & \longrightarrow & E_{n+1},\\
& X_n & \mapsto & D_i(X_n)=X_n\circ_i M
\end{array}$$
such that $D_0(X_0)=\mathbf{1}_E$ and for $X_n=(a_1\le\cdots\le
a_n)\in E_n,$ $X_n \circ_i M=(a_1\le\cdots\le a_i\le a_i+1\le
a_{i+1}+1\le\cdots\le a_n+1)$.
$$\begin{array}{cccc}
F_i :& E_n & \longrightarrow & E_{n-1} \\
& X_{n} & \mapsto & F_i(X)=X_{n}\circ_i X_0
\end{array}$$
such that for $X_n=(a_1\le\cdots\le a_n)\in E_n,$
$X\circ_iX_0=(a_1\le\cdots\le a_{i-1}\le a_{i+1}-1\le\cdots\le
a_n-1)$
\end{Pro}

\begin{proof}
To establish the above result, it is enough to show the following
properties.\\
$$\left\{ \begin{array}{lll}
F_iF_j=F_{j-1}F_i,\quad \mbox{if}\quad i\leq j\\
D_iD_j=D_{j+1}D_i,\quad \mbox{if} \quad i \leq j\\
F_iD_j=\left\{\begin{array}{lll}
D_{j-1}F_i,\quad \mbox{if}\quad i< j \\
id,\quad i\in\{j,j+1\}\\
D_jF_{i-1},\quad \mbox{if}\quad i> j+1
\end{array}\right.
\end{array}\right.$$
 From the associativity of the partial composition, one has  that for
$X_n\in E_n, Y_m\in E_m, Z_l\in E_l$ and for
 $i<j$ \\
$(X\circ_j Z)\circ_iY=(X\circ_iY)\circ_{j+m-1}Z,
(a)$.\\
Using relation $(a),$ we obtain:
\begin{enumerate}
\item[(1)] for $i< j$\\
$F_iF_j(X) = F_i(X\circ_jX_0) =(X\circ_jX_0)\circ_i X_0 = (X\circ_iX_0)\circ_{j-1}X_0=F_{j-1}F_i(X)$\\
$D_iD_j(X) = D_i(X\circ_j M) = (X\circ_j M)\circ_i M = (X\circ_i M)\circ_{j+1}M = D_{j+1} D_i(X)$\\
 $F_iD_j(X) = F_i(X\circ_j M)=(X\circ_j M)\circ_i X_0 = (X\circ_i X_0)\circ_{j-1} M = D_{j-1}F_i(X)$

\item[(2)] If $i>j+1,$ then $j<i-1.$ Using again relation (a), one has\\
$D_jF_{i-1}(X)=D_j(X\circ_{i-1}X_0)=(X\circ_{i-1}X_0)\circ_j M =
X\circ_j M)\circ_i X_0  =F_i D_j(x)$
\item[(3)] If $i\in\{j,j+1\}$, we have by direct calculation what
follows:
$$\begin{array}{ccc}
D_iD_i(X) & = &  (X\circ_i M)\circ_i M\\
 &= & a_1 \le \cdots \le
a_{i-1}\le a_i\le a_i+1\le a_{i+1} \\ & + & 1\le\cdots\le a_n+1)\circ_i M\\
& = & (a_1\le\cdots\le a_{i-1}\le a_i\le a_i\\  & + & 1\le a_i+2\le a_{i+1}+2\le\cdots\le a_n+2)\\
& = & D_{i+1}D_i(X)
\end{array}$$
$$\begin{array}{ccc}
F_i D_i(X) &= & (X\circ_iM)\circ_iX_0\\
&=&(a_1\le\cdots\le a_{i-1}\le a_i\le a_i+1\le a_{i+1}+1\le\cdots\le a_n+1)\circ_iX_0\\
 &=&(a_1\le\cdots\le a_{i-1}\le a_i+1-1\le a_{i+1}+1-1\le\cdots\le a_n+1-1)\\
 &=& X
\end{array}$$
$$\begin{array}{ccc}
F_{j+1}D_j(X)&=& (X\circ_jM)\circ_{j+1}X_0\\
&= & (a_1\le\cdots\le a_{j-1}\le a_j\le a_j+1\le a_{j+1}+1\le\cdots\le a_n+1)\circ_{j+1}X_0\\
& = & (a_1\le\cdots\le a_{j-1}\le a_j+1-1\le a_{j+1}+1-1\le\cdots\le
a_n+1-1)\\
& = & X
\end{array}$$
\end{enumerate}
\end{proof}
\begin{Rem}
Since the shift-operad
\end{Rem}

\end{document}